\documentclass[12pt,a4paper]{article}

\usepackage[english]{babel}
\usepackage{verbatim}
\usepackage{bm}

\usepackage{amsmath}
\usepackage{amssymb}
\usepackage{amsfonts}
\usepackage{amsbsy}

\usepackage{graphicx}
\usepackage{subfigure}
\usepackage{epsfig}

\usepackage{array}
\usepackage{booktabs}
\usepackage{caption}
\usepackage{tabularx}

\usepackage{morefloats}

\usepackage{color}

\begin{document}

\title{An accurate and efficient numerical framework \\
for adaptive numerical weather prediction}

\author{Giovanni Tumolo$^{(1)}$,\ \  Luca Bonaventura$^{(2)}$
 }

\maketitle

\begin{center}
{\small
$^{(1)}$
Earth System Physics Section \\
The Abdus Salam International Center for Theoretical Physics \\
Strada Costiera 11, 34151 Trieste, Italy \\
{\tt gtumolo@ictp.it}
}
\end{center}

\begin{center}
{\small
$^{(2)}$ MOX -- Modelling and Scientific Computing, \\
Dipartimento di Matematica ``F. Brioschi'', Politecnico di Milano \\
Via Bonardi 9, 20133 Milano, Italy\\
{\tt luca.bonaventura@polimi.it}
}
\end{center}

\date{}


\vspace*{0.5cm}


\vspace*{0.5cm}

\pagebreak

\abstract{We present an accurate and efficient discretization approach 
for the adaptive discretization of typical model equations employed in numerical weather prediction. 
A semi-Lagrangian approach is combined
with the TR-BDF2 semi-implicit time discretization method and with
a spatial discretization based on adaptive discontinuous finite elements.
The resulting method has full second order accuracy in time and can employ polynomial bases of arbitrarily high degree
 in space, is unconditionally stable and can effectively adapt the number of degrees
of freedom employed in each element, in order to balance accuracy and computational
cost. The $p-$adaptivity approach employed does not require remeshing, therefore it is  especially
suitable for applications, such as numerical weather prediction, in which a large number
of physical quantities are associated with a given mesh.
Furthermore, although the proposed method can be implemented on arbitrary unstructured and 
nonconforming meshes,
even its application on simple Cartesian meshes in spherical coordinates  can cure
effectively the pole problem by reducing the polynomial degree used in the polar elements.
 Numerical simulations of  classical  benchmarks for the shallow water and for the fully compressible Euler equations
  validate the method and  demonstrate its capability to achieve
 accurate results also at large Courant numbers, with time steps up to 100 times larger than
 those of typical explicit discretizations of the same problems, while reducing the computational
 cost thanks to the adaptivity algorithm.  
 }

\pagebreak

 \section{Introduction}
 \label{intro} \indent
 
 The  Discontinuous Galerkin (DG) spatial discretization approach is currently being 
 employed by an increasing number of   environmental fluid dynamics  models 
  see e.g. \cite{dawson:2006}, \cite{giraldo:2002}, \cite{lauter:2008},\cite{nair:2005},\cite{giraldo:2010},\cite{kelly:2012} and a more complete overview in \cite{bonaventura:2012}.
   This is motivated by the many attractive features of DG discretizations,
  such as high order accuracy, local mass  conservation and ease of massively parallel implementation.
  
 On the other hand, DG methods imply severe stability restrictions
 when coupled with explicit time discretizations.
 One traditional approach to overcome stability
 restrictions in low Mach number problems is the combination
  of  semi - implicit (SI) and semi - Lagrangian (SL) 
 techniques. In a series of papers \cite{restelli:2006}, 
  \cite{restelli:2009}, \cite{giraldo:2010}, \cite{dumbser:2013}, \cite{tumolo:2013} it has been shown that most of the computational 
  gains traditionally achieved in finite difference models by the application of SI, SL and SISL discretization methods
   are also attainable in the framework of DG approaches. 
  In particular,  in \cite{tumolo:2013} we have introduced a dynamically $p-$adaptive
 SISL-DG discretization approach for low Mach number 
  problems, that is quite effective in achieving high order spatial accuracy, while
 reducing substantially the computational cost.

In this paper, we apply the technique of  \cite{tumolo:2013}   to the shallow water equations
in spherical  geometry and to the  the fully compressible Euler equations, in order
to show its effectiveness for model problems typical of
global and regional weather forecasting.
 The advective form of the equations of motion is employed and the semi-implicit time discretization
is based on the TR-BDF2 method, see e.g. \cite{hosea:1996}, \cite{leveque:2007}.
This combination of two robust  ODE solvers yields a second order accurate, 
A-stable and L-stable method (see e.g. \cite{lambert:1991}), that is effective in damping selectively high frequency
modes. At the same time, it achieves full second order accuracy,
while the off-centering in the trapezoidal
rule, typically necessary  for realistic applications to nonlinear problems 
(see e.g. \cite{casulli:1994}, \cite{davies:2005}, \cite{tumolo:2013}), 
limits the accuracy in time to first order.
Numerical results presented in this paper show that the total computational
cost of one TR-BDF2  step is analogous to that of one step of the off-centered
trapezoidal rule, as well as the structure of the linear problems to be solved at each time step,
thus allowing to extend naturally to this more accurate method any implementation based on 
the off-centered trapezoidal rule. Numerical simulations of the  shallow water
   benchmarks proposed in \cite{williamson:1992}, \cite{lauter:2005}, \cite{jakob:1995} 
   and of the non-hydrostatic benchmarks proposed in \cite{skamarock:1994}, \cite{carpenter:1990}
 have been employed to validate the method and to demonstrate its capabilities.
 In particular, it will be shown that the present approach enables the use of time steps even 100
  times larger  than those allowed for DG models by standard explicit schemes, 
  see e.g. the results in \cite{nair:2005b}.

The method presented in this paper, just as its previous version in  \cite{tumolo:2013},
can be applied in principle on arbitrarily unstructured and even nonconforming meshes.
For example, a model based on this  method could run on a non conforming mesh of rectangular elements
built around the nodes of a reduced Gaussian grid \cite{hortal:1991}. For simplicity, however,
no such implementation has been developed so far. 
Here, only a simple Cartesian mesh has been used. If no degree adaptivity
is employed, this results in very high Courant numbers in the polar regions.
These do not  result in any special  stability problems for the present SISL discretization approach,
as it will be shown by the numerical results reported below.
On the other hand, even with an implementation based on a simple Cartesian
mesh in spherical coordinates, the flexibility of the DG space discretization allows to reduce the degree
of the basis and test functions employed close to the poles, thus making the effective
model resolution more uniform and solving the  efficiency issues related to the pole problem by static $p-$adaptivity.
This is especially advantageous because the  conditioning of the linear system to be solved at each time
step is greatly improved and, as a consequence, the number of iterations necessary for the
linear solver is reduced by approximately $80\%$, while at the same time
no spurious reflections nor artificial error increases are observed.

Beyond these computational advantages, we believe that the present approach based on
 $p-$adaptivity     is especially
suitable for applications to numerical weather prediction,  in contrast  to $h-$adaptivity  
approaches (that is, local mesh coarsening or refinement in which the size of some elements changes in time).
Indeed, in numerical weather prediction,  information that is necessary to carry out realistic simulations (such as  orography profiles, data on land use and soil type, land-sea masks) needs to be reconstructed on the computational mesh and has to be re-interpolated each time that
the mesh is changed. Furthermore, many physical parameterizations are highly sensitive to
the mesh size. Although devising better parameterizations that require less mesh-dependent tuning
is an important research goal,  more conventional
parameterizations will still be in use for quite some time. As a consequence, it is useful to improve the 
accuracy locally by adding supplementary degrees of freedom where necessary, as done in a  
 $p-$adaptive framework, without  having to change the underlying computational mesh.
In conclusion, the resulting modeling framework seems to be able to combine the
efficiency and high order accuracy of traditional SISL pseudo-spectral methods with the  
locality and flexibility of more standard DG approaches.  

In section \ref{shwater}, two examples of governing equations  are introduced.
 In section \ref{tr_rev}, the TR-BDF2 method is reviewed.
In section \ref{sphere} the approach employed for the advection of vector fields
in spherical geometry is described in detail.
In section  \ref{sisldg},  we introduce the SISL-DG discretization 
 approach for the shallow water equations in spherical geometry. 
 In section \ref{nhydro}, we outline its extension
 to the fully compressible Euler equations in a vertical plane.
 Numerical results are presented in section \ref{tests}, while in section \ref{conclu}
 we try to draw some conclusions and outline the path towards
  application of the concepts introduced here in the context of a non hydrostatic dynamical core.

\section{Governing equations}
\label{shwater} \indent

We consider as a basic model problem the two-dimensional shallow water equations 
on a rotating sphere (see e.g. \cite{gill:1982}).
These equations are a standard test bed for numerical
methods to be applied to the full equations of motion of atmospheric
or oceanic circulation models, see e.g. \cite{williamson:1992}. 
Among their possible solutions, they admit Rossby and inertial gravity waves,
as well as the response of such waves to orographic forcing.
We will use the  advective, vector form  of the shallow water equations:
\begin{eqnarray}
  && \frac{D h}{ Dt} = - h \nabla \cdot {\bf u}, \label{continuityeq}\\
  &&\frac{D {\bf u} }{ Dt} = - g \nabla h - f \hat{\bf k} \times {\bf u} -g \nabla b \label{vectmomentumeq}.
\end{eqnarray}
Here $h$ represents the fluid depth, $b$ the bathymetry elevation, 
$f$ the Coriolis parameter, $\hat{\bf k}$ the unit vector locally normal to the Earth's surface and 
$g$ the gravity force per unit mass on the Earth's surface. 
Assuming that $x,y$ are  orthogonal curvilinear coordinates
on the sphere (or on a portion of it), we denote by
 $m_x$ and $m_y$   the components of the (diagonal) metric tensor.
 Furthermore, we set
${\bf u}=(u,v)^T,$ where
$u$ and $v$ are the contravariant components of the velocity vector
in the coordinate direction $x$ and $y$ respectively,
multiplied by the corresponding metric tensor  components.
We also denote by
$\frac{D}{Dt}$ the Lagrangian derivative 
$$ \frac{D}{Dt} = \frac{\partial }{ \partial t} + 
   \frac{u}{m_x}  \frac{\partial }{ \partial x} +
   \frac{v}{m_y}  \frac{\partial }{ \partial y} ,$$
so that $u = m_x \frac{D x}{Dt}, v = m_y \frac{D y}{Dt}$. 
In particular, in this paper standard spherical coordinates will be employed.

 As an example of a more complete model, we will also consider
the fully compressible, non hydrostatic equations of motion.
Following e.g.  \cite{cullen:1990},  \cite{bonaventura:2000}, \cite{davies:2005}, they can be written as
 \begin{eqnarray*}
   && \frac{D \Pi}{Dt} = - \left( \frac{c_p}{c_v}-1 \right) \Pi \nabla \cdot {\bf u}, \\
   && \frac{D {\bf u} }{Dt} = - c_p \Theta \nabla \Pi -g \hat{\bf k}, \\
   && \frac{D \Theta}{Dt} = 0.
  \end{eqnarray*} 
 where, being $p_0 $ a reference pressure value, $\Theta = T \big( \frac{p}{p_0} \big)^{-R/c_p}$ is the potential temperature, 
 $ \Pi = \big( \frac{p}{p_0} \big)^{R/c_p}$ is the Exner pressure,  while
 $c_p, c_v, R $ are the constant pressure and constant volume specific heats and the gas constant of dry air
  respectively. Here the Coriolis force is omitted for simplicity. Notice also that, by a slight abuse of notation,
  in the three-dimensional case  ${\bf u}=(u,v,w)^T $ denotes the three dimensional velocity
  field and the $\frac{D}{Dt},$  $\nabla $ operators are also three-dimensional, while we will assume ${\bf u}=(u,w)^T $ in the 
  description of $(x,z)$ two dimensional, vertical slice models.
   It is customary to rewrite such equations in terms of perturbations with respect to a steady
  hydrostatic reference profile, so that assuming
  $ \Pi(x,y,z,t) = \pi^*(z) + \pi(x,y,z,t), $ $ \Theta(x,y,z,t) = \theta^*(z) + \theta(x,y,z,t) $ with
 $ \hspace{2.5mm} c_p \theta^* \frac{d \pi^*}{d z} = -g, $ one obtains for a vertical plane
  \begin{eqnarray}
   && \frac{D \Pi}{Dt} = - \left( \frac{c_p}{c_v}-1 \right) \Pi \nabla \cdot { \bf u}, \label{vslice_conteq}\\
   && \frac{D u }{Dt}  = - c_p \Theta \frac{\partial \pi}{\partial x}, \label{vslice_ueq}\\
   && \frac{D w }{Dt}  = - c_p \Theta \frac{\partial \pi}{\partial z} + g \frac{\theta}{\theta^*}, \label{vslice_veq} \\
   && \frac{D \theta}{Dt} =  - \frac{d \theta^*}{dz} w. \label{vslice_eneq}
  \end{eqnarray} 
 It can be observed that equations \eqref{vslice_conteq}-\eqref{vslice_eneq} are isomorphic
 to equations \eqref{continuityeq}-\eqref{vectmomentumeq}, which will allow to extend
 almost automatically the discretization approach proposed for the former to the more general
 model.

 \section{Review of the TR-BDF2 method} 
  \label{tr_rev}
 
\begin{figure}[hcb]
\begin{center}
 \includegraphics[height=0.35\textheight]{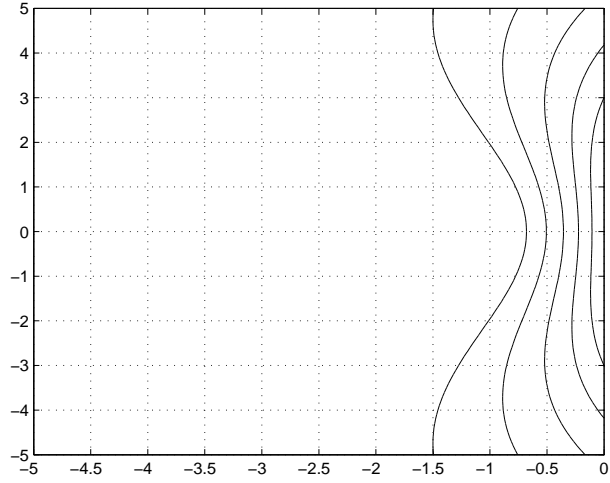} 
 \end{center}
 \caption{Contour levels of the absolute value of the stability
 function of the TR-BDF2 method without off centering in the first stage.
 Contour spacing is $0.1$ from $0.5$ to $1.$}
 \label{stabreg_trbdf_noff}
\end{figure}

\begin{figure}[t]
\begin{center}
 \includegraphics[height=0.35\textheight]{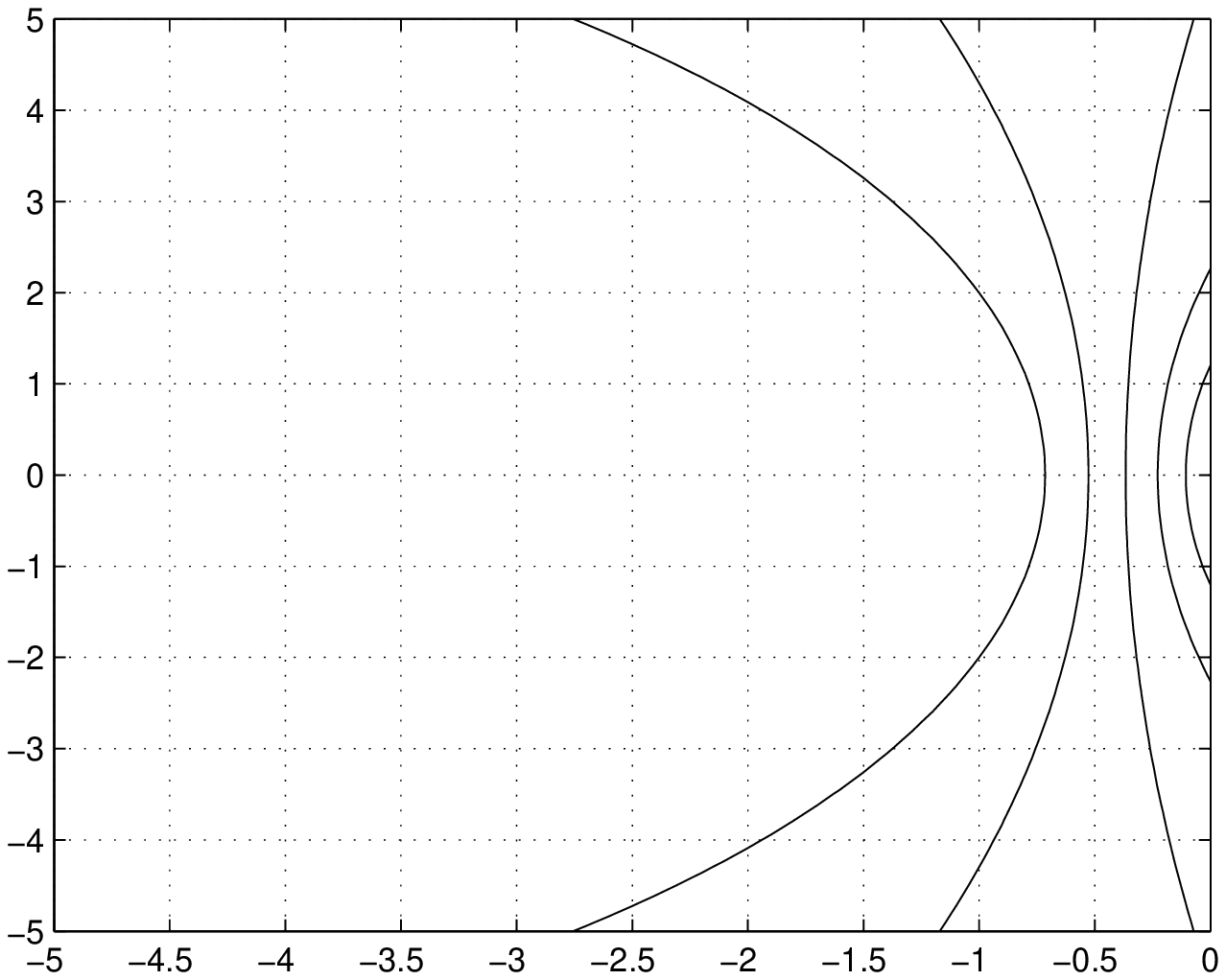} 
  \end{center}
 \caption{Contour levels of the absolute value of the stability
 function of the off centered Crank-Nicolson method with averaging parameter
 $\theta=0.6 $ (equivalent to an off centering parameter valued $0.05$).  Contour spacing is $0.1$ from $0.5$ to $1.$}
 \label{stabreg_theta06}
\end{figure}

\begin{figure}[b]
\begin{center}
 \includegraphics[height=0.35\textheight]{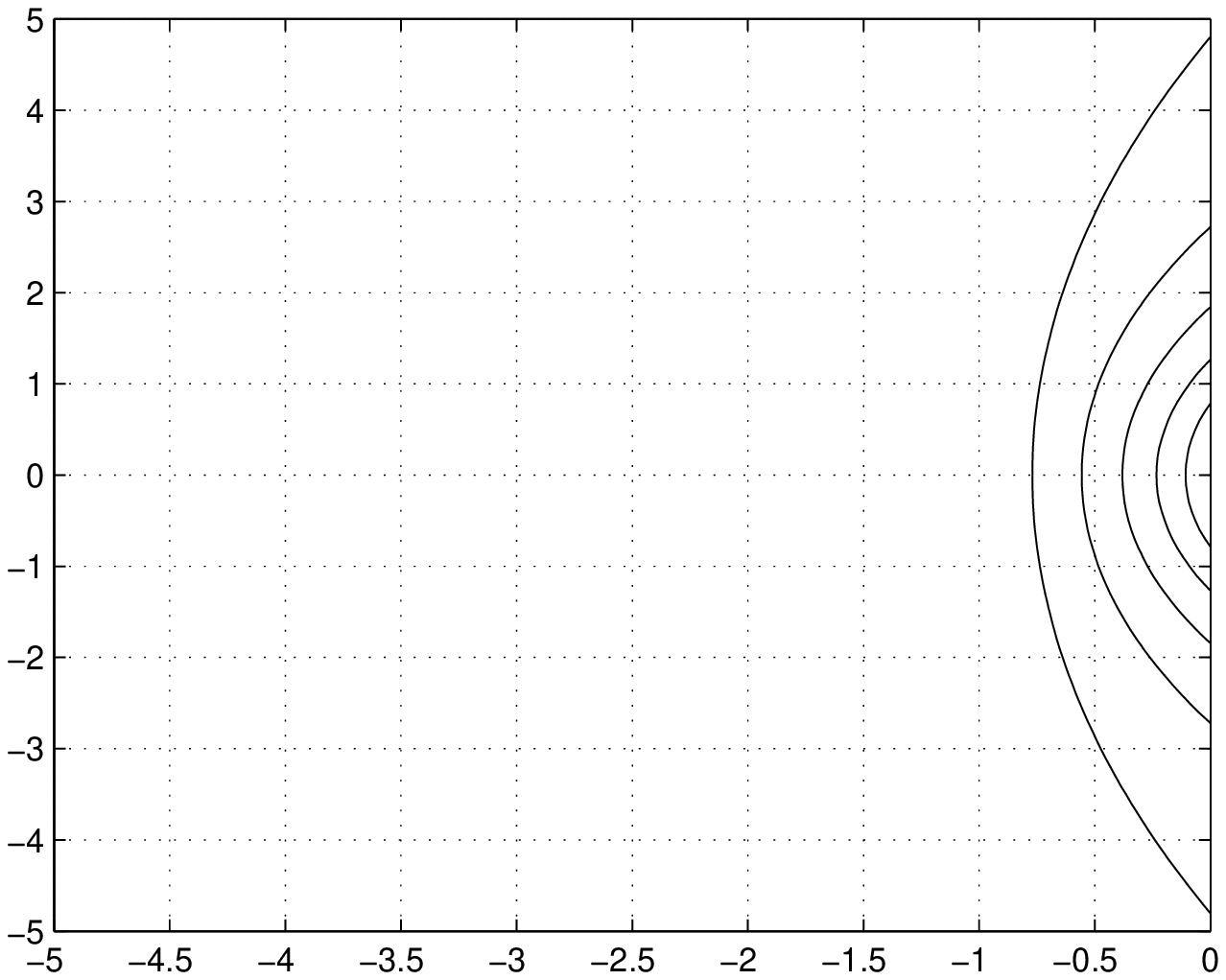} 
  \end{center}
 \caption{Contour levels of the absolute value of the stability
 function of the off centered Crank-Nicolson method with averaging parameter
 $\theta=0.7 $  (equivalent to an off centering parameter valued $0.1$). Contour spacing is $0.1$ from $0.5$ to $1.$}
 \label{stabreg_theta07}
\end{figure}

\begin{figure}[t]
\begin{center}
 \includegraphics[height=0.33\textheight]{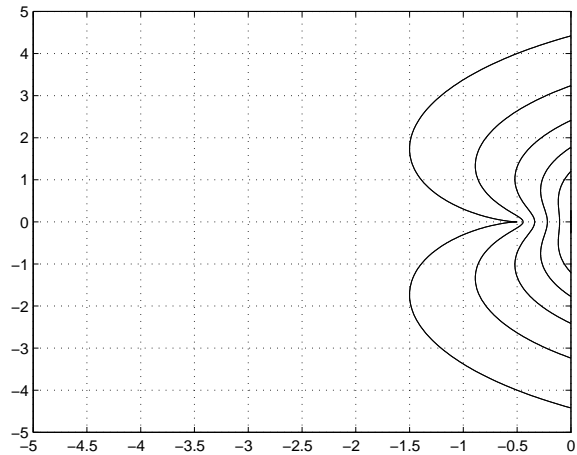} 
   \end{center}
 \caption{Contour levels of the absolute value of the stability
 function of the BDF2 method.
 Contour spacing is $0.1$ from $0.5$ to $1.$}
 \label{bdf_sreg}
\end{figure}

\begin{figure}[b]
\begin{center}
 \includegraphics[height=0.35\textheight]{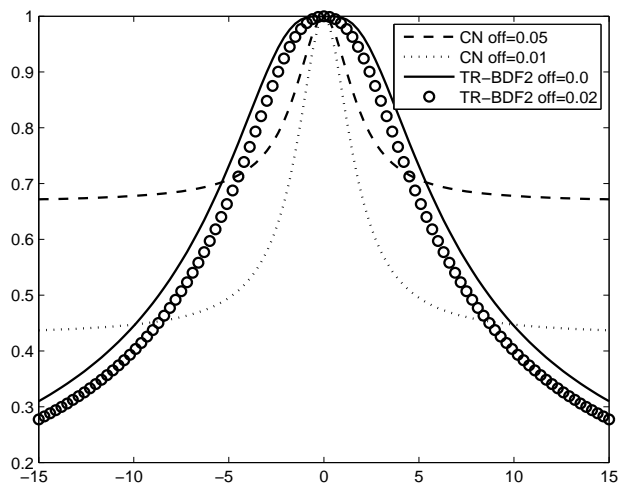} 
   \end{center}
 \caption{Graph of the absolute value of the stability
 functions of several L-stable methods along the imaginary axis.}
 \label{section_imaxis}
\end{figure}

 We review here some properties of the so called
 TR-BDF2 method, which was first introduced in  \cite{bank:1985}.  
  Given a Cauchy problem
 \begin{eqnarray}
 \label{cauchypb}
   {\bf y}^{\prime}&=&{\bf f}({\bf y},t)    \nonumber \\
   {\bf  y}(0) &=& {\bf y}_0
\end{eqnarray}
   and considering  a time 
 discretization employing  a constant time step $\Delta t,$ the TR-BDF2 method is defined by
 the two following implicit stages:
\begin{eqnarray}\label{trbdf2}
  {\bf u}^{n+2\gamma} - \gamma \Delta t  {\bf f}({\bf u}^{n+2\gamma},t_{n}+2\gamma\Delta t)
   &=& {\bf u}^n + \gamma \Delta t {\bf f}({\bf u}^{n},t_{n}), \nonumber \\
  {\bf u}^{n+1} - \gamma_2 \Delta t  {\bf f}({\bf u}^{n+1},t_{n+1}) &=& (1-\gamma_3 ){\bf u}^n +\gamma_3 {\bf u}^{n+2\gamma}. 
\end{eqnarray}
Here $\gamma \in [0,1/2] $ is an implicitness parameter and $$\gamma_2 =  \frac{1-2\gamma}{2(1-\gamma)}, 
 \ \ \ \gamma_3  =\frac{1-\gamma_2}{2\gamma} .$$
 It is immediate that the first stage  of \eqref{trbdf2} is simply the application of the trapezoidal rule
 (or Crank-Nicolson method) over the interval 
 $[t_n,t_{n}+2\gamma\Delta t].$
 It could also be substituted by an off centered  Crank-Nicolson step
 without reducing the overall accuracy of the method.
 The outcome of this stage is then used to turn the two step BDF2 method into a single step, two stages method.
 This combination of two robust  stiff solvers yields a method with several interesting accuracy and stability
 properties, that were   analyzed in detail in \cite{hosea:1996}. As shown in this paper,
 this analysis is most easily carried out by rewriting the method as

\begin{eqnarray}
{\bf k}_1& =&  {\bf f}\left ({\bf u}^n, t_{n}\right) \nonumber \\
 {\bf k}_2 &=& {\bf f}\left ({\bf u}^n +\gamma \Delta t {\bf k}_1  +\gamma \Delta t {\bf k}_2,   t_{n}+\gamma \Delta t\right)\nonumber \\
 {\bf k}_3 &=& {\bf f}\left ({\bf u}^n +\frac{1-\gamma}2\Delta t{\bf k}_1+\frac{1-\gamma}2\Delta t{\bf k}_2 
 +\gamma\Delta t{\bf k}_3,t_{n+1} \right )
 \nonumber \\ 
 {\bf u}^{n+1} &=& {\bf u}^n +\Delta t \left (\frac{1-\gamma}2{\bf k}_1+\frac{1-\gamma}2{\bf k}_2 +\gamma{\bf k}_3\right).
\label{tr_dirk}
\end{eqnarray}
In this formulation, the TR-BDF2 method is clearly a Singly Diagonal Implicit Runge Kutta (SDIRK) method, so that one can rely on the theory for this class  of methods to derive
stability and accuracy results (see e.g. \cite{lambert:1991}). 
Notice that the same method has been rediscovered in \cite{butcher:2000} and has been 
analyzed and applied also in \cite{giraldo:2013}, to treat the implicit terms in the
framework of an Additive Runge Kutta approach (see e.g. \cite{kennedy:2003}).
As shown in  \cite{hosea:1996},
the TR-BDF2 method is second order accurate and A-stable  for any value of $\gamma.$ 
Written as in \eqref{tr_dirk}, the method can also be proven to constitute a (2,3) embedded Runge-Kutta
pair, with companion coefficients given by
$$
(1-\frac{\sqrt{2}}4)/3, \ \ \ (1+3\frac{\sqrt{2}}4)/3,  \ \ \  \frac{2-\sqrt{2}}6,
$$
provided that no off centering is employed in the first stage of \eqref{trbdf2}.
This equips the method with an extremely efficient estimator of the time discretization error.
Furthermore, for $\gamma=1-\sqrt{2}/2$ it is also  L-stable. Therefore, with this coefficient
value it can be safely applied to 
 problems with eigenvalues whose imaginary part is large, such as typically arise from the discretization
 of  hyperbolic problems.  This is not the case for the standard 
trapezoidal rule (or Crank-Nicolson) implicit method, whose linear stability region is exactly bounded by
the imaginary axis. As a consequence, it is common to apply the trapezoidal rule with
off centering, see e.g. \cite{casulli:1994}, \cite{davies:2005} as well as \cite{tumolo:2013},  which results in a first
order time discretization. TR-BDF2 appears therefore to be an interesting one step alternative to maintain
full second order accuracy, especially considering that, if formulated as \eqref{trbdf2}, it is equivalent
to performing two Crank-Nicolson steps with slightly modified coefficients.
In order to highlight the advantages of the proposed method in terms of accuracy with respect
to other common robust stiff solvers, we plot in figure \ref{stabreg_trbdf_noff} the
contour levels of the absolute value of the   linear stability
 function of the TR-BDF2 method without off centering in the first stage, compared to the analogous
 contours of the off centered Crank-Nicolson method with averaging parameter
 $\theta=0.6, \theta=0.7 $ in figures \ref{stabreg_theta06},  \ref{stabreg_theta07}, respectively,
 and to those of the BDF2 method in figure \ref{bdf_sreg}. It is immediate to see that
 TR-BDF2 introduces less damping around the imaginary axis for moderate values
 of the time step. On the other hand, TR-BDF2 is more selective in damping very large eigenvalues,
  as clearly displayed in figure \ref{section_imaxis},
 where the absolute values of the  linear  stability
 functions of the same methods (with the exception of BDF2, for which an
 explicit representation of the stability function is not available) are plotted along the imaginary axis.

 \pagebreak

\pagebreak

\clearpage
 
  \section{Review of Semi-Lagrangian evolution operators for vector fields on the sphere} 
  \label{sphere}
  
  The semi-Lagrangian method can be described introducing the concept of 
  evolution operator, along the lines of \cite{morton:1995,morton:1998}.
  Indeed, let $G = G({\bf x}, t)$ denote a generic function of space and time
  that is the solution of
  $$ \frac{D G}{Dt} = \frac{\partial G}{ \partial t} + 
   \frac{u}{m_x}  \frac{\partial G}{ \partial x} +
   \frac{v}{m_y}  \frac{\partial G}{ \partial y} =0.$$
To approximate this solution on the time interval
 $[t^n,t^{n+1}],$ a numerical evolution operator $E$ is introduced,
    that approximates the exact evolution operator associated to 
the frozen velocity field ${\bf u}^{*}=(u^*,v^*)^T ,$ that may coincide
with the velocity field at time level $t^n$ or with an extrapolation derived
from more   previous time levels.
More precisely, if ${\bf X}(t;t^{n+1},{\bf x})$ denotes the solution of 
\begin{equation}
\frac{d {\bf X}(t;t^{n+1},{\bf x} )}{dt}={\bf u}^*({\bf X}(t;t^{n+1},{\bf x})) 
\label{lagode}
\end{equation}
with initial datum ${\bf X}(t^{n+1};t^{n+1},{\bf x})={\bf x}$ at time $t=t^{n+1}$,
then the expression $ [E(t^{n},\Delta t) G]({\bf x}) $ denotes a numerical approximation of $G^n({\bf x}_D)$
where ${\bf x}_D={\bf X}(t^{n};t^{n+1},{\bf x})$
and the notation $G^n({\bf x}) = G({\bf x}, t^n)$ is used.
Since ${\bf x}_D$ is
nothing but the position at time $t^n$ of the fluid parcel reaching location ${\bf x}$ at time $t^{n+1}$,
according to standard terminology, it is called the departure point associated with the arrival point ${\bf x}$.
Different  methods can be employed to approximate ${\bf x}_D$; 
in this paper, for simplicity, the method proposed in \cite{mcgregor:1993} has been employed in spherical geometry.
Furthermore, to guarantee an accuracy compatible with that of the semi-implicit time-discretization,
 an extrapolation ${\bf u}^{n+\frac 12}$ of the velocity field at the intermediate time level $t^n+\Delta t/2 $ 
 was used as ${\bf u}^*$ in (\ref{lagode}). On the other hand, in the application to   Cartesian geometry (for the vertical slice discretization), 
 a simple first order Euler method with sub-stepping was employed, see e.g. \cite{giraldo:1999}, \cite{rosatti:2005}.

In case of the advection of a vector field
$$ \frac{D {\bf G}}{Dt} = \frac{\partial {\bf G}}{ \partial t} + 
   \frac{u}{m_x}  \frac{\partial {\bf G}}{ \partial x} +
   \frac{v}{m_y}  \frac{\partial {\bf G}}{ \partial y} =0,$$
as in the momentum equation  (\ref{vectmomentumeq}),
the extension of this approach has to take into account the curvature of the spherical manifold.
More specifically, unit basis vectors at the departure point are not in general aligned  
with those at the arrival point, i.e., if $\hat{\boldsymbol{i}},\hat{\boldsymbol{j}},\hat{\boldsymbol{k}}$
represent a unit vector triad, in general 
$ \hspace{1mm} \hat{\bf i}({\bf x}) \neq \hat{\boldsymbol{i}}({\bf x}_D),  
  \hspace{1mm} \hat{\boldsymbol{j}}({\bf x}) \neq \hat{\boldsymbol{j}}({\bf x}_D), 
  \hspace{1mm} \hat{\boldsymbol{k}}({\bf x}) \neq \hat{\boldsymbol{k}}({\bf x}_D). $
  
To deal with this issue two approaches are available.
The first, intrinsically Eulerian, consists in the introduction of the Christoffel symbols in the 
covariant derivatives definition, giving rise to the well known metric terms, before the SISL discretization, 
and then in the approximation along the trajectories of those metric terms.
This approach has been shown to be source of instabilities in a semi-Lagrangian frame, see e.g.
\cite{ ritchie:1988, cote:1988, cote:1988b, desharnais:1990} and therefore is not adopted in this work.
The second approach, more suitable for semi-Lagrangian discretizations, takes into account  the curvature of the manifold
  only at discrete level, i.e.  after the SISL discretization has been performed.
Many variations of this idea have been proposed, see e.g. \cite{ ritchie:1988, cote:1988, cote:1988b, bates:1990, temperton:2001}.
In \cite{staniforth:2010}, they have all been derived in a unified way by the introduction of a proper rotation matrix 
that transforms vector components in the departure-point unit vector triad 
$\hat{\boldsymbol{i}}_D=\hat{\boldsymbol{i}}({\bf x}_D),$
$\hat{\boldsymbol{j}}_D=\hat{\boldsymbol{j}}({\bf x}_D),$
$\hat{\boldsymbol{k}}_D=\hat{\boldsymbol{k}}({\bf x}_D)$
into vector components in the arrival-point unit vector triad
$\hat{\boldsymbol{i}}=\hat{\boldsymbol{i}}({\bf x}),$
$\hat{\boldsymbol{j}}=\hat{\boldsymbol{j}}({\bf x}),$
$\hat{\boldsymbol{k}}=\hat{\boldsymbol{k}}({\bf x})$.
To see how this rotation matrix comes into play, it is sufficient to consider the action of the evolution operator $E$ 
on a given vector valued function of space and time $\boldsymbol{G}$, defined as an approximation of 
\begin{equation}
  \left[ E(t^{n}, \Delta t) \boldsymbol{G} \right] ({\bf x}) = \boldsymbol{G}^n({\bf x}_D), 
  \label{evoloponvectors}
\end{equation}
and to write this equation componentwise.
$\boldsymbol{G}^n({\bf x}_D)$ is known through its components in the departure point unit vector triad:
\begin{equation}
  \boldsymbol{G}^n({\bf x}_D) = \mathcal{G}^n_x({\bf x}_D)\hat{\boldsymbol{i}}_D+
                                       \mathcal{G}^n_y({\bf x}_D)\hat{\boldsymbol{j}}_D+
                                       \mathcal{G}^n_z({\bf x}_D)\hat{\boldsymbol{k}}_D.
  \label{vectGexpasion}                                     
\end{equation}                                        
                                        
Therefore, via (\ref{evoloponvectors}), the components of 
$\left[ E(t^{n}, \Delta t) \boldsymbol{G} \right] ({\bf x})$
in the unit vector triad at the same point are given by projection of (\ref{vectGexpasion}) along
$ \hat{\boldsymbol{i}},$ $\hat{\boldsymbol{j}},$ $ \hat{\boldsymbol{k}}$:

\begin{equation}
  \hspace{1mm}\hat{\boldsymbol{i}} \cdot \boldsymbol{G}^n({\bf x}_D) =
  \mathcal{G}^n_x({\bf x}_D)\hspace{1mm}\hat{\boldsymbol{i}} \cdot\hat{\boldsymbol{i}}_D+
  \mathcal{G}^n_y({\bf x}_D)\hspace{1mm}\hat{\boldsymbol{i}}  \cdot\hat{\boldsymbol{j}}_D+
  \mathcal{G}^n_z({\bf x}_D)\hspace{1mm}\hat{\boldsymbol{i}} \cdot\hat{\boldsymbol{k}}_D, \nonumber 
\end{equation}
\begin{equation}
  \hspace{1mm}\hat{\boldsymbol{j}}  \cdot \boldsymbol{G}^n({\bf x}_D) =
  \mathcal{G}^n_x({\bf x}_D)\hspace{1mm}\hat{\boldsymbol{j}}  \cdot\hat{\boldsymbol{i}}_D+
  \mathcal{G}^n_y({\bf x}_D)\hspace{1mm}\hat{\boldsymbol{j}}  \cdot\hat{\boldsymbol{j}}_D+
  \mathcal{G}^n_z({\bf x}_D)\hspace{1mm}\hat{\boldsymbol{j}} \cdot\hat{\boldsymbol{k}}_D, \nonumber
\end{equation}
\begin{equation}
  \hspace{1mm}
  \hat{\boldsymbol{k}} \cdot \boldsymbol{G}^n({\bf x}_D) =
  \mathcal{G}^n_x({\bf x}_D)\hat{\boldsymbol{k}} \cdot\hat{\boldsymbol{i}}_D+
  \mathcal{G}^n_y({\bf x}_D)\hat{\boldsymbol{k}}  \cdot\hat{\boldsymbol{j}}_D+
  \mathcal{G}^n_z({\bf x}_D)\hat{\boldsymbol{k}}  \cdot\hat{\boldsymbol{k}}_D, \nonumber
\end{equation}
i.e., in matrix notation:
\begin{equation*}
  \begin{pmatrix}
    \hspace{1mm}\hat{\boldsymbol{i}}  \cdot \left[ E(t^{n}, \Delta t) \boldsymbol{G} \right]({\bf x}) \\
    \hspace{1mm}\hat{\boldsymbol{j}} \cdot \left[ E(t^{n}, \Delta t) \boldsymbol{G} \right]({\bf x}) \\
    \hspace{1mm}\hat{\boldsymbol{k}}\cdot \left[ E(t^{n}, \Delta t) \boldsymbol{G} \right]({\bf x})
  \end{pmatrix}
  = {\bf R} \begin{pmatrix}  
    \mathcal{G}^n_x \\
    \mathcal{G}^n_y \\
    \mathcal{G}^n_z
  \end{pmatrix} \ \ \ \ \ \ \ {\rm where }
  \end{equation*}
 \begin{equation*}
 {\bf R} =   \begin{bmatrix}
                 \hat{\boldsymbol{i}}  \cdot\hat{\boldsymbol{i}}_D
                  & \hat{\boldsymbol{i}}  \cdot\hat{\boldsymbol{j}}_D&
                   \hat{\boldsymbol{i}} \cdot\hat{\boldsymbol{k}}_D \\
                 \hat{\boldsymbol{j}}  \cdot\hat{\boldsymbol{i}}_D & 
                 \hat{\boldsymbol{j}}\cdot\hat{\boldsymbol{j}}_D & \hat{\boldsymbol{j}}  \cdot\hat{\boldsymbol{k}}_D \\
                 \hat{\boldsymbol{k}}  \cdot\hat{\boldsymbol{i}}_D & \hat{\boldsymbol{k}} \cdot\hat{\boldsymbol{j}}_D & \hat{\boldsymbol{k}}  \cdot\hat{\boldsymbol{k}}_D \\
               \end{bmatrix}.
\end{equation*}
Under the shallow atmosphere approximation \cite{thuburn:2013}, $\bf R$ 
can be reduced to the $2\times2 $ rotation matrix
\begin{equation}
{\boldsymbol \Lambda } = {\boldsymbol \Lambda }({\bf x},{\bf x}_D) =  
   \begin{bmatrix}
                \Lambda_{11} &  \Lambda_{12}  \\  
                \Lambda_{21} &  \Lambda_{22} 
   \end{bmatrix}, 
   \label{lambdadef}
\end{equation}
where, as shown in \cite{staniforth:2010}, $\Lambda_{11}=\Lambda_{22}=(R_{11}+R_{22})/(1+R_{33}), 
\hspace{1mm} \Lambda_{12}=-\Lambda_{21}=(R_{12}-R_{21})/(1+R_{33}).$
Therefore, in the following the evolution operator for vector fields will be defined componentwise as
 
\begin{equation}
  \begin{pmatrix}
    \hspace{1mm}\hat{\boldsymbol{i}}  \cdot \left[ E(t^{n}, \Delta t) \boldsymbol{G} \right]({\bf x}) \\
    \hspace{1mm}\hat{\boldsymbol{j}} \cdot \left[ E(t^{n}, \Delta t) \boldsymbol{G} \right]({\bf x})  
  \end{pmatrix}
  =   {\boldsymbol \Lambda } \begin{pmatrix}  
    \mathcal{G}^n_x({\bf x}_D) \\
    \mathcal{G}^n_y ({\bf x}_D)
  \end{pmatrix}.
  \label{evoloponvectors2}
\end{equation}

 \section{A novel SISL time integration approach  for the shallow water equations on the sphere} 
  \label{sisldg}

 The  SISL discretization of equations.(\ref{continuityeq})-(\ref{vectmomentumeq}) based on (\ref{trbdf2})
 is then obtained by performing the two stages in (\ref{trbdf2}) after reinterpretation of the intermediate
 values in a semi-Lagrangian fashion. Furthermore, in order to avoid the solution of a nonlinear system,
 the dependency on $h$ in $h \nabla \cdot {\bf u}$ is linearized in time, as common
 in semi-implicit discretizations based on the trapezoidal rule, see e.g. \cite{casulli:1994},\cite{tumolo:2013}.
 Numerical experiments reported in the following show that this does not prevent to achieve
 second order accuracy in the regimes of interest for numerical weather prediction.
 The TR stage of the SISL time semi-discretization of the equations in vector form 
 (\ref{continuityeq})-(\ref{vectmomentumeq}) is given by
\begin{eqnarray}\label{TR_SWEcontinuityeq}
   h^{n+2\gamma} &+& \gamma \Delta t \hspace{1.0mm} h^n \hspace{1.0mm} \nabla \cdot {\bf u}^{n+2\gamma}\nonumber \\
   & = &  
   E\left(t^n, 2 \gamma \Delta t\right) \left[ h -\gamma \Delta t \hspace{1.0mm} h \hspace{1.0mm} \nabla \cdot {\bf u} \right], 
\end{eqnarray}
 
\begin{eqnarray}\label{TR_SWEmomentumeq}
   && {\bf u}^{n+2\gamma}  
   + \gamma \Delta t \Big[  g \nabla h^{n+2\gamma} + f \hat{\bf k} \times { \bf u}^{n+2\gamma}\Big] 
   = -\gamma \Delta t \hspace{1.0mm} g\nabla b \nonumber \\
   && + E\big(t^n, 2 \gamma \Delta t\big) \left\{ {\bf u} -\gamma \Delta t \left[ g ( \nabla h + \nabla b ) + f \hat{\bf k} \times {\bf u} \right] \right\}. 
\end{eqnarray} 
 The TR stage is then followed by the BDF2 stage:

\begin{eqnarray}\label{BDF2_SWEcontinuityeq}
   h^{n+1} &+& \gamma_2 \Delta t \hspace{1.0mm} h^{n+2\gamma} \hspace{1.0mm} \nabla \cdot {\bf u}^{n+1} \nonumber  \\
   &= & \big( 1-\gamma_3\big) E\big(t^n,\Delta t\big) h \nonumber \\
   &+& \gamma_3 E\big(t^n+2\gamma\Delta t,(1-2\gamma)\Delta t\big) h,  
\end{eqnarray}

 \begin{eqnarray}\label{BDF2_SWEmomentumeq}
   &&{\bf u}^{n+1}  
   + \gamma_2 \Delta t \Big[  g \nabla h^{n+1} + f \hat{\bf k} \times {\bf u}^{n+1}\Big] 
   = - \gamma_2 \Delta t \hspace{1.0mm} g \nabla b \nonumber \\
   && + \big( 1-\gamma_3\big) E\big(t^n,\Delta t\big)  {\bf u} + 
         \gamma_3 E\big(t^n+2\gamma\Delta t,(1-2\gamma)\Delta t\big) {\bf u}.  
\end{eqnarray} 
 For each of the two stages, the spatial discretization can be performed along the lines described
 in \cite{tumolo:2013}, allowing for variable polynomial order to locally represent the solution in each element.
 The spatial discretization approach considered
is independent of the nature of the mesh and could also
be implemented for fully unstructured and even non conforming meshes.
For simplicity, however, in this paper  only an implementation
on a structured mesh in longitude-latitude coordinates has been developed.
  In principle,
 either Lagrangian or hierarchical Legendre bases could be employed. 
We will work almost exclusively with hierarchical bases, because they provide a natural environment
for the implementation of a $p-$adaptation algorithm, see for example 
\cite{zienkiewicz:1983}. A central issue in  finite element formulations for fluid problems is the choice 
of appropriate approximation spaces for the velocity and pressure variables (in the 
context of SWE, the role of the  pressure is played by the free surface elevation). 
An inconsistent choice of the two approximation spaces indeed may result in a solution
 that is polluted by spurious modes, for the specific case of  SWE
see for example \cite{leroux:2005,walters:1983a,walters:1983b} as well as the more recent 
and comprehensive analysis in \cite{leroux:2013}.
Here, we have not investigated this issue in depth, but the model implementation allows
for approximations of higher polynomial degree $ p^u $ for the velocity fields than $ p^h $ for the height field.
Even though no systematic study was performed, no significant differences were
noticed between results obtained with equal or unequal degrees. In the following,
only results with unequal degrees $  p^u =p^h+1 $ are reported, with the exception
of an empirical convergence test for a steady geostrophic flow.

 All the integrals appearing in the elemental equations are evaluated
by means of Gaussian numerical quadrature formulae, with a number of
quadrature nodes consistent with the local polynomial degree being used. 
In particular, notice that integrals of terms in the image of the evolution operator $E,$ i.e. 
of functions evaluated at the departure points of the trajectories arriving at the quadrature nodes, 
cannot be computed exactly (see e.g. \cite{morton:1988,priestley:1994}), since such functions are not polynomials. 
Therefore a sufficiently accurate approximation of these integrals is needed, 
which may entail the need to employ numerical quadrature formulae with more nodes than
the minimal requirement implied by the local polynomial degree. This overhead 
 is actually compensated by the fact that, for each Gauss node, the computation of the departure point 
 is only to be executed once for all the quantities to be interpolated.

 After spatial discretization has been performed, the discrete degrees of freedom representing velocity unknowns
 can be replaced in the respective discrete height equations, yielding in each case
 a linear system whose structure is entirely analogous to that obtained in \cite{tumolo:2013}.
 The  non-symmetric linear systems obtained from the TR-BDF2 stages are solved in our implementation
by  the GMRES method \cite{saad:1986}.  A classical stopping criterion
based on a relative error tolerance  of $10^{-10} $ was employed (see e.g. \cite{kelley:1995}).
For the GMRES solver, so far, only a block
diagonal preconditioning was employed. As it will be shown in section \ref{tests}, the condition number of the systems to be solved
can be greatly reduced if lower degree elements are employed close to the poles.
In any case, the total computational
cost of one TR-BDF2  step is entirely analogous to that of one step of the standard off centered
trapezoidal rule employed in \cite{tumolo:2013},  since the structure of the systems is the same
but for each stage only a fraction of the time step is being computed.
 Once ${h}^{n+1}$ has been computed by solving this linear system,  then $ {\bf u}^{n+1}$
  can be recovered by back substituting into the momentum equation.

\section{Extension of the time integration approach to the Euler equations} 
\label{nhydro}

  In this section, we show that the previously proposed method can be extended
  seamlessly to the fully compressible Euler equations as formulated in equations 
  (\ref{vslice_conteq}) - (\ref{vslice_eneq}). For simplicity, only the application to the
  $(x,z)$ two dimensional vertical slice case is presented, but the extension to
  three dimensions is straightforward.
  Again, in order to avoid the solution of a nonlinear system,
 the dependency on $\Pi$ in $\Pi \nabla \cdot {\bf u}$ and the dependency on $\Theta$ in $\Theta \nabla \pi$ are
 linearized in time, as common in semi-implicit discretizations based on the trapezoidal rule, see e.g. \cite{cullen:1990}, \cite{bonaventura:2000}.
  
 The semi-Lagrangian counterpart of the TR substep of (\ref{trbdf2}) is first applied to
to (\ref{vslice_conteq}) - (\ref{vslice_eneq}), so as to obtain:

\begin{eqnarray}\label{TR_VSLcontinuityeq}
   && \pi^{n+2\gamma} 
   + 
   \gamma \Delta t \hspace{1.0mm} \left( c_p/c_v - 1 \right) \Pi^n \nabla \cdot {\bf u}^{n+2\gamma} = -\pi^* \nonumber \\
   && + E\left(t^n, 2 \gamma \Delta t\right) \left[ \Pi -\gamma \Delta t \left( c_p/c_v - 1 \right) \Pi \hspace{1.0mm}  \nabla \cdot {\bf u} \right], 
\end{eqnarray}
 
\begin{eqnarray}\label{TR_VSLmomentumeq_x}
   && u^{n+2\gamma}  
   + \gamma \Delta t \hspace{1.0mm} c_p \Theta^n \frac{\partial \pi}{\partial x}^{n+2\gamma} = \nonumber \\
   && E(t^n, 2\gamma \Delta t) \left[ u - \gamma \Delta t \hspace{1.0mm} c_p \Theta \frac{\partial \pi}{\partial x} \right], 
\end{eqnarray} 

\begin{eqnarray}\label{TR_VSLmomentumeq_z}
   &&w^{n+2\gamma}  
   + \gamma \Delta t \left( c_p \Theta^n \frac{\partial \pi}{\partial z}^{n+2\gamma} 
    - g \frac{\theta^{n+2\gamma} }{ \theta^* } \right) = \nonumber \\
   && E(t^n, 2\gamma \Delta t) \left[ w - \gamma \Delta t \left( c_p \Theta \frac{\partial \pi}{\partial z} - g \frac{\theta}{ \theta^* }\right) \right], 
\end{eqnarray} 

\begin{equation}\label{TR_VSLenereq}
   \theta^{n+2\gamma}  
   + \gamma \Delta t \frac{d \theta^*}{d z} w^{n+2\gamma} 
   = E(t^n, 2\gamma \Delta t) \left[ \theta - \gamma \Delta t \frac{d \theta^*}{d z} w \right].
\end{equation} 

Following \cite{cullen:1990} the time semi-discrete energy equation (\ref{TR_VSLenereq}) 
can be inserted into the time semi-discrete vertical momentum equation (\ref{TR_VSLmomentumeq_z}), 
in order to decouple the momentum and the energy equations as follows

\begin{eqnarray}\label{TR_VSLmomentumeq_z+ener}
&& \left( 1 + (\gamma \Delta t)^2 \frac{g}{\theta^*} \frac{d \theta^*}{d z} \right) w^{n+2\gamma}  
+ \gamma \Delta t c_p \Theta^n \frac{\partial \pi}{\partial z}^{n+2\gamma} = \nonumber \\
&&  E(t^n, 2\gamma \Delta t) \left[ w - \gamma \Delta t \left( c_p \Theta \frac{\partial \pi}{\partial z} - g \frac{\theta}{ \theta^* }\right) \right]  \nonumber \\ 
&& + \gamma \Delta t \frac{g}{\theta^*} E(t^n, 2\gamma \Delta t) \left[ \theta - \gamma \Delta t \frac{d \theta^*}{d z} w \right]. 
\end{eqnarray} 
Equations (\ref{TR_VSLcontinuityeq}), (\ref{TR_VSLmomentumeq_x}) and (\ref{TR_VSLmomentumeq_z+ener})
are a set of three equations in three unknowns only, namely $\pi, u,$ and $w$  that can be compared
with equations (\ref{TR_SWEcontinuityeq}), (\ref{TR_SWEmomentumeq}) with $f=0$  
and $m_x=m_y=1$ (Cartesian geometry).
From the comparison it is clear that the two formulations are isomorphic under
 correspondence $\pi \longleftrightarrow h, u \longleftrightarrow u, w \longleftrightarrow v.$
 
We can then consider the semi-Lagrangian counterpart of the 
BDF2 substep of (\ref{trbdf2}) applied to (\ref{vslice_conteq}) - (\ref{vslice_eneq})
to obtain:

\begin{eqnarray}\label{BDF2_VSLcontinuityeq}
   \pi^{n+1} 
   &+& 
   {\gamma}_2 \Delta t \left( c_p/c_v - 1 \right) \Pi^{n+2\gamma} \nabla \cdot {\bf u}^{n+1} \nonumber  \\
   &=& -\pi^* + (1 - {\gamma}_3) [ E\left(t^n, \Delta t\right)  \Pi ] \nonumber \\
   &+& {\gamma}_3 [ E\left(t^n+2\gamma \Delta t, (1-2\gamma)\Delta t\right)  \Pi ], 
\end{eqnarray}
 
\begin{eqnarray}\label{BDF2_VSLmomentumeq_x}
   u^{n+1}  
   &+& {\gamma}_2 \Delta t \hspace{1.0mm} c_p \Theta^{n+2\gamma} \frac{\partial \pi}{\partial x}^{n+1} \nonumber \\
   &=& (1 - {\gamma}_3) [ E\left(t^n, \Delta t\right)  u ] \nonumber \\
   &+& {\gamma}_3 [ E\left(t^n+2\gamma \Delta t, (1-2\gamma)\Delta t\right)  u ],  
\end{eqnarray} 

\begin{eqnarray}\label{BDF2_VSLmomentumeq_z}
   w^{n+1}  
   &+& {\gamma}_2 \Delta t \left( c_p \Theta^{n+2\gamma} \frac{\partial \pi}{\partial z}^{n+1} 
    - g \frac{\theta^{n+1} }{ \theta^* } \right)\nonumber \\
   &=& (1 - {\gamma}_3) [ E\left(t^n, \Delta t\right) w ] \nonumber \\
   &+& {\gamma}_3 [ E\left(t^n+2\gamma \Delta t, (1-2\gamma)\Delta t\right) w ],  
\end{eqnarray} 

\begin{eqnarray}\label{BDF2_VSLenereq}
   \theta^{n+1}  
  &+& {\gamma}_2 \Delta t \frac{d \theta^*}{d z} w^{n+1} \nonumber\\ 
  &=& (1 - {\gamma}_3) [ E\left(t^n, \Delta t\right) \theta ] \nonumber \\
  &+& {\gamma}_3 [ E\left(t^n+2\gamma \Delta t, (1-2\gamma)\Delta t\right) \theta ]. 
\end{eqnarray} 
Again, following \cite{cullen:1990}, the time semi-discrete energy equation (\ref{BDF2_VSLenereq}) 
can be inserted into the time semi-discrete vertical momentum equation (\ref{BDF2_VSLmomentumeq_z}), in order to decouple the momentum and the energy equations:

\begin{eqnarray}\label{BDF2_VSLmomentumeq_z+ener}
  && \left( 1 + ({\gamma}_2 \Delta t)^2 \frac{g}{\theta^*} \frac{d \theta^*}{d z} \right) w^{n+1}  
     + {\gamma}_2 \Delta t  \hspace{1.0mm} c_p \Theta^{n+2\gamma} \frac{\partial \pi}{\partial z}^{n+1} = \nonumber \\
  && (1 - {\gamma}_3) [ E\left(t^n, \Delta t\right) w ] + {\gamma}_3 [ E\left(t^n+2\gamma \Delta t, (1-2\gamma)\Delta t\right) w ] +  \\
  && {\gamma}_2 \Delta t \frac{g}{\theta^*} \left\{ (1 - {\gamma}_3) [ E\left(t^n, \Delta t\right) \theta ] + {\gamma}_3 [ E\left(t^n+2\gamma \Delta t, (1-2\gamma)\Delta t\right) \theta ] \right\}. \nonumber
\end{eqnarray} 

Now equations (\ref{BDF2_VSLcontinuityeq}), (\ref{BDF2_VSLmomentumeq_x}) and (\ref{BDF2_VSLmomentumeq_z+ener})
are a set of three equations in three unknowns only, namely $\pi, u,$ and $w,$  that can be compared
with equations (\ref{BDF2_SWEcontinuityeq}), (\ref{BDF2_SWEmomentumeq}) with $f=0$  
and $m_x=m_y=1$ (Cartesian geometry).  Again, it is easy to see that also in this case exactly the same structure results
 as in equations (\ref{BDF2_SWEcontinuityeq})-(\ref{BDF2_SWEmomentumeq}) with the correspondence
$\pi \longleftrightarrow h, u \longleftrightarrow u, w \longleftrightarrow v$,
so that the  approach (and code) proposed for the shallow water equations can be extended to the fully compressible Euler equation in a straightforward way.
 
\clearpage 
 
\section{Numerical experiments}
\label{tests}

The numerical method introduced in section \ref{sisldg}
has been implemented and tested on a number of relevant test cases
using different initial conditions and bathymetry profiles,
in order to assess its accuracy and stability properties and 
to analyze the impact of the  $p-$adaptivity strategy.
Whenever a reference solution was available,
the relative errors were computed in the $L^1,L^2 $ and $L^\infty $ norms 
at the final time $t_f$ of the simulation according to \cite{williamson:1992} as: 

\begin{eqnarray}\label{errornorms}
&& l_1(h) = \frac{I \left[ \hspace{1mm} \left| h(\cdot, t_f) - h_{ref}(\cdot, t_f) \right| \hspace{1mm} \right]}
                 {I \left[ \hspace{1mm} \left| h_{ref}(\cdot, t_f) \right| \hspace{1mm} \right] }, \\ 
&& l_2(h) = \frac{ \Bigl\{ I \Bigl[ \hspace{1mm} \bigl( \hspace{1mm} h(\cdot, t_f) - h_{ref}(\cdot, t_f) \hspace{1mm} \bigr)^2 \hspace{1mm} \Bigr] \Bigr\}^{1/2}}
                 { \bigl\{ I \bigl[ \hspace{1mm} h_{ref}(\cdot, t_f)^2 \hspace{1mm} \bigr] \bigr\}^{1/2} }, \\
&& l_{\infty}(h) = \frac{\max \hspace{1mm} \left| h(\cdot, t_f) - h_{ref}(\cdot, t_f) \right| }
                        {\max \hspace{1mm} \left| h_{ref}(\cdot, t_f) \right| },
\end{eqnarray}
where $h_{ref}$ denotes the reference solution for a model variable $h$ and $I$ is a discrete approximation of the global integral 
\begin{equation}
 I(h)= \frac{\int_{\Omega} \, h \, m_x m_y \, d{\bf x}}{\int_{\Omega} \, m_x m_y \,  d{\bf x}},
 \label{normalizedintegral}
\end{equation}
computed by an appropriate numerical quadrature rule, consistent
with the numerical approximation being tested, and the  maximum is computed over all nodal values.

The test cases considered for the shallow water equations 
in spherical geometry are
\begin{itemize}
\item a steady-state geostrophic flow: in particular, we have analyzed results in test case 2 of  \cite{williamson:1992}  
 in the configuration least favorable for methods employing longitude-latitude meshes;
\item the unsteady flow with exact analytical solution described in \cite{lauter:2005}; 
\item the polar rotating low-high, introduced in \cite{mcdonald:1989},   aimed at
showing that no problems arise even in the case of strong cross polar flows;
\item zonal flow over an isolated mountain and Rossby-Haurwitz wave of wavenumber 4,
corresponding respectively to test cases 5 and 6 in \cite{williamson:1992}.
\end{itemize}
For the first two tests,
analytic solutions are available and empirical convergence tests can be performed.
The test cases considered for the discretization of equations \eqref{vslice_conteq}-\eqref{vslice_eneq} are
\begin{itemize}
\item inertia gravity waves involving the evolution of a potential temperature
      perturbation in a channel with periodic boundary conditions and 
      uniformly stratified environment with constant Brunt-W\"ais\"al\"a frequency, 
      as described in \cite{skamarock:1994}; 
\item a rising thermal bubble given by the evolution of a warm bubble 
      in a constant potential temperature environment, as described in \cite{carpenter:1990}.
\end{itemize}

In all the numerical experiments performed for this paper, neither spectral filtering 
nor explicit diffusion of any kind were employed, the only numerical diffusion
being implicit in the time discretization approach.
We have not yet investigated to which extent the quality of the solutions
is affected by this choice, but 
this should be taken into account when comparing quantitatively the results of the present method
to those of reference models, such as the one described in \cite{jakob:1995}, in which
explicit numerical diffusion is added. Sensitivity of the comparison results
to the amount of numerical diffusion has been highlighted in several model validation exercises,
 see e.g. \cite{ripodas:2009}.
 
Since semi-implicit, semi-Lagrangian methods are most efficient for low  Froude 
number flows, where the typical
velocity is much smaller than that of the fastest propagating 
waves, all the tests considered fall in this hydrodynamical regime. Therefore,
in order to assess the method efficiency, a distinction has been made between 
the maximum Courant number based on the velocity, on one hand, and, on the other hand,
the maximum Courant number based on the celerity, or
the maximum Courant number based on the sound speed,
defined respectively as
$$  C_{vel}= \max \frac{\|{\bf u}\|_{\infty}\Delta t}{\Delta x / p} $$
$$    C_{cel}= \max \frac{\sqrt{g h}\Delta t}{\Delta x / p}, \ \ 
    C_{snd}= \max \frac{\sqrt{(c_p/c_v) R \Theta \Pi }\Delta t}{\Delta x / p},$$
    where $\Delta x$ is to be interpreted as generic value of the meshsize in either coordinate
    direction.
For the tests in which $p-$adaptivity was employed,
if $p_I^n$ denotes the local polynomial degree used at timestep $t^n $ to represent a model variable
inside the $I-th$ element of the mesh, while $p_{max}$ is the 
maximum local polynomial degree considered, 
the efficiency of the method in reducing the computational effort has been measured by
monitoring the evolution of the quantities
$$\Delta_{dof}^n = \frac{ \sum_{I=1}^N (p_I^n +1)^2 }{ N (p_{max}+1)^2 },
 \ \ \ \ \ \Delta_{iter}^n = \frac{{\rm ITN}^n_{adapt}}{{\rm ITN}^n_{max}},$$
where $N$ is the total number of elements,
  ${\rm ITN}^n_{adapt} $ denotes the total number
of GMRES iterations at   time step $n $ for the adapted local degrees configuration
and ${\rm ITN}^n_{max}$  the total number
of GMRES iterations at   time step $n $ for the configuration with maximum  degree in all elements,
respectively. Average values of these indicators over the simulations performed are reported in the
following,  denoted by $\Delta_{dof}^{average}$ and $\Delta_{iter}^{average}$ respectively.
The error between the adaptive solution and the corresponding one obtained 
with uniform maximum polynomial degree everywhere has been measured in terms of (\ref{errornorms}).
Finally, in some cases conservation of global invariants has been monitored by evaluating
 at each time step the following 
global integral quantities:

\begin{equation}
 J(q^n) = \frac{ I(q(\cdot,t^n)) - I(q(\cdot,t^0)) }{I(q(\cdot,t^0))},
\end{equation}

where $I(q)$ has been defined in (\ref{normalizedintegral})
and $q^n=q(\cdot,t^n)$ is the density associated to each global invariant.
According to the choice of $q$, following invariants are considered: 
mass, i.e. $q=q_{mass}=h$, 
total energy, i.e. $q=q_{energ}=\frac{1}{2} ( h \boldsymbol{u} \cdot \boldsymbol{u} + g (h^2 - b^2)$),
and potential enstrophy, i.e. $q=q_{enstr}=\frac{1}{2h} (\hat{\boldsymbol{k}} \cdot \nabla \times \boldsymbol{u}+f)^2.$ 

 \subsection{Steady-state geostrophic flow}
\label{test2}
We first  consider the test case 2 of \cite{williamson:1992}, 
where the solution is a steady state flow
with velocity field corresponding to a zonal solid body rotation and 
$h$ field obtained from the velocity ones through geostrophic balance.
 All the parameter values are taken as in \cite{williamson:1992}. 
The flow orientation parameter has been chosen here as $\alpha = \pi/2 -0.05,$ 
making the test more challenging on a longitute-latitude mesh. 
Error norms associated to the solution obtained on a mesh of $10 \times 5$ elements for different polynomial degrees are
shown in tables \ref{t2convrate_h_tab}, \ref{t2convrate_u_tab} and \ref{t2convrate_v_tab} for $h, $ $u $ and $v,$ respectively.
All the results have been computed at $t_f = 10 $ days at fixed maximum Courant numbers $C_{cel}=8, C_{vel}=2, $ so that
different values of  $\Delta t $ have been employed for different polynomial order.
We remark that the resulting time steps are significantly larger than those
allowed by typical explicit time discretizations for analogous DG space discretizations,
see e.g. the results in \cite{nair:2005b}. 
The spectral decay in the error norms can be clearly observed, until the time error becomes dominant.
For  better comparison with the results in \cite{nair:2005b}, we consider again the configuration with $p^h=6, p^u=7$ on $10 \times 5$ elements,
which corresponds to the same resolution in space as for  the $150 \times 8 \times 8$ grid used in \cite{nair:2005b}.
While  $\Delta t = 36 $ s is used in \cite{nair:2005b} giving a $l_{\infty}(h) \approx 8 \times 10^{-6},$ 
the proposed SISLDG formulation can be run with $\Delta t = 3600 $ s, in which case $l_{\infty}(h) \approx 3 \times 10^{-7},$ 
and the average number of iterations required by the linear solver is 1 for the TR substep and 4 for the BDF2 substep.

 \begin{table}[htbc]
\[
\begin{array}{cccccc}
\toprule
p^h & p^u & \Delta t \ [s] &  l_1(h)                & l_2(h)                 & l_{\infty}(h)        \\
\midrule
2   & 3   &  4800          & 5.558 \times 10^{-3}   & 6.805 \times 10^{-3}   & 1.914 \times 10^{-2} \\
3   & 4   &  3600          & 6.017 \times 10^{-4}   & 8.176 \times 10^{-4}   & 2.569 \times 10^{-3} \\
4   & 5   &  2880          & 1.743 \times 10^{-5}   & 2.405 \times 10^{-5}   & 9.024 \times 10^{-5} \\
5   & 6   &  2400          & 1.586 \times 10^{-6}   & 2.281 \times 10^{-6}   & 1.058 \times 10^{-5} \\
6   & 7   &  2057          & 8.829 \times 10^{-8}   & 1.206 \times 10^{-7}   & 4.926 \times 10^{-7} \\
7   & 8   &  1800          & 1.246 \times 10^{-8}   & 1.590 \times 10^{-8}   & 4.158 \times 10^{-8} \\
8   & 9   &  1600          & 5.641 \times 10^{-9}   & 5.952 \times 10^{-9}   & 6.320 \times 10^{-9} \\
\bottomrule
\end{array}
\]
\caption{Relative errors on $h $ for different polynomial degrees,   SWE test case 2 with $\alpha = \pi/2 -0.05 $
at time $t_f = 10 $ days. }
\label{t2convrate_h_tab}
\end{table}

\begin{table}[htbc]
\[
\begin{array}{cccccc}
\toprule
p^h & p^u & \Delta t \ [s] &   l_1(u)               & l_2(u)                 & l_{\infty}(u)        \\
\midrule
2   & 3   &  4800          & 6.351 \times 10^{-2}   & 6.432 \times 10^{-2}   & 1.143 \times 10^{-1} \\
3   & 4   &  3600          & 9.505 \times 10^{-3}   & 1.037 \times 10^{-2}   & 2.106 \times 10^{-2} \\
4   & 5   &  2880          & 4.288 \times 10^{-4}   & 4.887 \times 10^{-4}   & 2.393 \times 10^{-3} \\
5   & 6   &  2400          & 4.598 \times 10^{-5}   & 4.830 \times 10^{-5}   & 1.706 \times 10^{-4} \\
6   & 7   &  2057          & 2.057 \times 10^{-6}   & 2.262 \times 10^{-6}   & 5.879 \times 10^{-6} \\
7   & 8   &  1800          & 2.162 \times 10^{-7}   & 2.358 \times 10^{-7}   & 6.428 \times 10^{-7} \\
8   & 9   &  1600          & 2.013 \times 10^{-8}   & 2.276 \times 10^{-8}   & 3.268 \times 10^{-8} \\
\bottomrule
\end{array}
\]
\caption{Relative errors on $u $ for different polynomial degrees,   SWE test case 2 with $\alpha = \pi/2 -0.05 $
at time $t_f = 10 $ days.}
\label{t2convrate_u_tab}
\end{table}

\begin{table}[htbc]
\[
\begin{array}{cccccc}
\toprule
p^h & p^u & \Delta t \ [s] &   l_1(v)               & l_2(v)                 & l_{\infty}(v)        \\
\midrule
2   & 3   &  4800          & 1.001 \times 10^{-1}   & 1.016 \times 10^{-1}   & 2.698 \times 10^{-1} \\
3   & 4   &  3600          & 1.859 \times 10^{-2}   & 1.823 \times 10^{-2}   & 6.848 \times 10^{-2} \\
4   & 5   &  2880          & 7.376 \times 10^{-4}   & 7.428 \times 10^{-4}   & 2.884 \times 10^{-3} \\
5   & 6   &  2400          & 8.185 \times 10^{-5}   & 8.307 \times 10^{-5}   & 2.574 \times 10^{-4} \\
6   & 7   &  2057          & 3.074 \times 10^{-6}   & 3.173 \times 10^{-6}   & 1.123 \times 10^{-5} \\
7   & 8   &  1800          & 3.370 \times 10^{-7}   & 3.432 \times 10^{-7}   & 1.323 \times 10^{-6} \\
8   & 9   &  1600          & 2.175 \times 10^{-8}   & 2.317 \times 10^{-8}   & 5.124 \times 10^{-8} \\
\bottomrule
\end{array}
\]
\caption{Relative errors on $v $ for different polynomial degrees,   SWE test case 2 with $\alpha = \pi/2 -0.05 $
at time $t_f = 10 $ days.}
\label{t2convrate_v_tab}
\end{table}
Another convergence test was performed for $p^h = p^u = 3, $  increasing the number of elements and correspondingly
decreasing the value of the time step. In this case, the maximum Courant numbers
vary because of the mesh inhomogeneity, so that $ 2 < C_{cel} < 18,$  $ 0.5 < C_{vel} < 4.$
 The results are reported in tables 
\ref{t2convrate_h_tab_pfix}, \ref{t2convrate_u_tab_pfix} and \ref{t2convrate_v_tab_pfix} for $h, $ $u $ and $v,$ respectively. The empirical convergence order $q_2^{emp}$ based on the $l_2 $ norm errors has also been estimated, showing that in this stationary test convergence rates above the second order of the time discretization can be achieved.

\begin{table}[htbc]
\[
\begin{array}{ccccccc}
\toprule
N_x  \times N_y  & \Delta t \ [s]    &  l_1(h)                & l_2(h)                 & l_{\infty}(h)   & q_2^{emp}     \\
\midrule
10   \times  5   &  3600            & 2.557 \times 10^{-4}   & 3.495 \times 10^{-4}   & 1.403 \times 10^{-3} & - \\
20   \times 10   &  1800            & 2.187 \times 10^{-5}   & 2.889 \times 10^{-5}   & 1.566 \times 10^{-4} & 3.6 \\
40   \times 20   &   900            & 2.530 \times 10^{-6}   & 3.353 \times 10^{-6}   & 1.430 \times 10^{-5} &  3.1  \\
80   \times 40   &   450            & 3.996 \times 10^{-7}   & 5.534 \times 10^{-7}   & 3.134 \times 10^{-6} & 2.6 \\
\bottomrule
\end{array}
\]
\caption{Relative errors on $h$ for different number of elements, $p^h = p^u = 3, $   SWE test case 2 with $\alpha = \pi/2 -0.05 $
at time $t_f = 10 $ days.}
\label{t2convrate_h_tab_pfix}
\end{table}

\begin{table}[htbc]
\[
\begin{array}{ccccccc}
\toprule
N_x \times N_y  & \Delta t \ [s] &   l_1(u)               & l_2(u)                 & l_{\infty}(u)    & q_2^{emp}     \\
\midrule
10  \times  5   &  3600       & 2.769 \times 10^{-3}   & 3.358 \times 10^{-3}   & 8.948 \times 10^{-3} &  - \\
20  \times 10   &  1800       & 2.896 \times 10^{-4}   & 3.720 \times 10^{-4}   & 2.414 \times 10^{-3} & 3.2 \\
40  \times 20   &   900       & 3.647 \times 10^{-5}   & 4.563 \times 10^{-5}   & 2.473 \times 10^{-4} & 3.0 \\
80  \times 40   &   450       & 6.826 \times 10^{-6}   & 1.035 \times 10^{-5}   & 9.525 \times 10^{-5} & 2.1 \\
\bottomrule
\end{array}
\]
\caption{Relative errors on $u $ for different number of elements, $p^h = p^u = 3, $   SWE test case 2 with $\alpha = \pi/2 -0.05 $
at time $t_f = 10 $ days.}
\label{t2convrate_u_tab_pfix}
\end{table}

\begin{table}[htbc]
\[
\begin{array}{cccccc}
\toprule
N_x \times N_y & \Delta t \ [s] &   l_1(v)               & l_2(v)                 & l_{\infty}(v)  & q_2^{emp}       \\
\midrule
10  \times  5   &  3600       & 3.309 \times 10^{-3}   & 3.346 \times 10^{-3}   & 8.250 \times 10^{-3} &  - \\
20  \times 10   &  1800       & 4.016 \times 10^{-4}   & 4.233 \times 10^{-4}   & 1.255 \times 10^{-3} & 3.0  \\
40  \times 20   &   900       & 5.180 \times 10^{-5}   & 5.578 \times 10^{-5}   & 2.329 \times 10^{-4} & 2.9 \\
80  \times 40   &   450       & 9.405 \times 10^{-6}   & 1.214 \times 10^{-5}   & 7.763 \times 10^{-5} & 2.2 \\
\bottomrule
\end{array}
\]
\caption{Relative errors on $v $ for different number of elements, $p^h = p^u = 3, $ SWE test case 2 with $\alpha = \pi/2 -0.05 $
at time $t_f = 10 $ days.}
\label{t2convrate_v_tab_pfix}
\end{table}

\clearpage

\subsection{Unsteady flow with analytic solution}
\label{lauter}

In a second, time dependent test, the analytic solution of
(\ref{continuityeq})-(\ref{vectmomentumeq})  derived in \cite{lauter:2005}  has been employed to assess
the performance of the proposed discretization. More specifically, the analytic solution defined in
formula  (23)  of \cite{lauter:2005}  was used. Since  the exact solution is periodic, 
the initial profiles also correspond to the exact solution an integer number of days later.
The proposed SISLDG scheme has been integrated up to $t_f= 5 $ days with
$p^h =4 $ and $ p^u = 5 $  on meshes with increasing number of elements,
while the time step has been decreased accordingly.
In this case, the maximum Courant numbers vary because of the mesh dishomogeneity, so that $ 4 < C_{cel} < 26,$  $ 1.25 < C_{vel} < 8.$
Error norms for $h, u, v$ of the above-mentioned integrations have been computed at $t_f=5 $ days and 
displayed in tables \ref{tlauterconvrate_h_tab} - \ref{tlauterconvrate_v_tab}.  An empirical order estimation
shows that full second order accuracy in time is attained.
 
\begin{table}[htbc]
\begin{center}
\begin{tabular}{cccccc}
  \toprule
  $N_x \times N_y$          &  $\Delta t  \ [s]$ & $l_1(h)$               &  $l_2(h)$               & $l_{\infty}(h)$        & $q_2^{emp} $     \\
  \midrule
  $10 \times \hspace{2mm}5$ &    3600            & $5.456 \times 10^{-3}$ &  $6.120 \times 10^{-3}$ & $9.537 \times 10^{-3}$ & - \\
  $20 \times 10           $ &    1800            & $1.246 \times 10^{-3}$ &  $1.397 \times 10^{-3}$ & $2.143 \times 10^{-3}$ & 2.1 \\
  $40 \times 20           $ &    900             & $3.039 \times 10^{-4}$ &  $3.410 \times 10^{-4}$ & $5.207 \times 10^{-4}$ &2.0 \\
  $80 \times 40           $ &    450             & $7.548 \times 10^{-5}$ &  $8.475 \times 10^{-5}$ & $1.292 \times 10^{-4}$ & 2.0\\  
  \bottomrule
\end{tabular}
\caption{Relative errors on $h $ at different resolutions,   L\"auter test case.}
\label{tlauterconvrate_h_tab}
\end{center}
\end{table}

\begin{table}[htbc]
\begin{center}
\begin{tabular}{cccccc}
  \toprule
  $N_x \times N_y$          & $\Delta t  \ [s]$ & $l_1(u)$               & $l_2(u)$               & $l_{\infty}(u)$        & $q_2^{emp}$     \\
  \midrule
  $10 \times \hspace{2mm}5$ &    3600           & $6.567 \times 10^{-2}$ & $7.848 \times 10^{-2}$ & $1.670 \times 10^{-1}$ & - \\      
  $20 \times 10           $ &    1800           & $1.665 \times 10^{-2}$ & $1.994 \times 10^{-2}$ & $3.931 \times 10^{-2}$ &2.0 \\
  $40 \times 20           $ &    900            & $4.210 \times 10^{-3}$ & $5.032 \times 10^{-3}$ & $9.811 \times 10^{-3}$ & 2.0 \\
  $80 \times 40           $ &    450            & $1.057 \times 10^{-3}$ & $1.261 \times 10^{-3}$ & $2.452 \times 10^{-3}$ & 2.0 \\  
  \bottomrule
\end{tabular}
\caption{Relative errors on $u $ at different resolutions,   L\"auter test case.}
\label{tlauterconvrate_u_tab}
\end{center}
\end{table}

\begin{table}[htbc]
\begin{center}
\begin{tabular}{cccccc}
  \toprule
  $N_x \times N_y$          & $\Delta t  \ [s]$ & $l_1(v)$               & $l_2(v)$               & $l_{\infty}(v)$        & $q_2^{emp}$     \\
  \midrule
  $10 \times \hspace{2mm}5$ &    3600           & $1.174 \times 10^{-1}$ & $1.198 \times 10^{-1}$ & $2.316 \times 10^{-1}$ & -\\      
  $20 \times 10           $ &    1800           & $2.939 \times 10^{-2}$ & $3.002 \times 10^{-2}$ & $5.561 \times 10^{-2}$ & 2.0\\
  $40 \times 20           $ &    900            & $7.336 \times 10^{-3}$ & $7.497 \times 10^{-3}$ & $1.390 \times 10^{-2}$ &2.0\\
  $80 \times 40           $ &    450            & $1.833 \times 10^{-3}$ & $1.874 \times 10^{-3}$ & $3.464 \times 10^{-3}$ & 2.0 \\  
  \bottomrule
\end{tabular}
\caption{Relative errors on $v $ at different resolutions,   L\"auter test case.}
\label{tlauterconvrate_v_tab}
\end{center}
\end{table}
For comparison, analogous errors have been computed with the same discretization parameters
but employing the off centered Crank Nicolson method of \cite{tumolo:2013} with $\theta=0.6$.
 The resulting improvement in the errors between the TRBDF2 scheme and the off-centered Crank Nicolson is 
 achieved at an essentially equivalent computational cost in terms of total CPU time employed.
   
 \begin{table}[htbc]
\begin{center}
\begin{tabular}{cccccc}
  \toprule
  $N_x \times N_y$          &  $\Delta t  \ [s]$ & $l_1(h)$               &  $l_2(h)$               & $l_{\infty}(h)$        & $q_2^{emp} $     \\
  \midrule
  $10 \times \hspace{2mm}5$ &    3600            & $1.444 \times 10^{-2}$ &  $1.633 \times 10^{-2}$ & $2.398 \times 10^{-2}$ & -   \\
  $20 \times 10           $ &    1800            & $8.742 \times 10^{-3}$ &  $9.894 \times 10^{-3}$ & $1.445 \times 10^{-2}$ & 0.7 \\
  $40 \times 20           $ &    900             & $4.814 \times 10^{-3}$ &  $5.451 \times 10^{-3}$ & $7.956 \times 10^{-3}$ & 0.9 \\
  $80 \times 40           $ &    450             & $2.526 \times 10^{-3}$ &  $2.861 \times 10^{-3}$ & $4.177 \times 10^{-3}$ & 0.9 \\  
  \bottomrule
\end{tabular}
\caption{Relative errors on $h $ at different resolutions,   L\"auter test case with off centered Crank Nicolson, $\theta=0.6$.}
\label{tlauterconvrate_h_tab_thet06}
\end{center}
\end{table}

\begin{table}[htbc]
\begin{center}
\begin{tabular}{cccccc}
  \toprule
  $N_x \times N_y$          & $\Delta t  \ [s]$ & $l_1(u)$               & $l_2(u)$               & $l_{\infty}(u)$        & $q_2^{emp}$     \\
  \midrule
  $10 \times \hspace{2mm}5$ &    3600           & $1.800 \times 10^{-1}$ & $2.092 \times 10^{-1}$ & $3.810 \times 10^{-1}$ & - \\      
  $20 \times 10           $ &    1800           & $1.077 \times 10^{-1}$ & $1.255 \times 10^{-1}$ & $2.155 \times 10^{-1}$ & 0.7 \\
  $40 \times 20           $ &    900            & $5.895 \times 10^{-2}$ & $6.880 \times 10^{-2}$ & $1.186 \times 10^{-1}$ & 0.9 \\
  $80 \times 40           $ &    450            & $3.084 \times 10^{-2}$ & $3.603 \times 10^{-2}$ & $6.234 \times 10^{-2}$ & 0.9 \\  
  \bottomrule
\end{tabular}
\end{center}
\caption{Relative errors on $u $ at different resolutions,   L\"auter test case with off centered Crank Nicolson, $\theta=0.6$.}
\label{tlauterconvrate_u_tab_thet06}
\end{table}

\begin{table}[htbc]
\begin{center}
\begin{tabular}{cccccc}
  \toprule
  $N_x \times N_y$          & $\Delta t  \ [s]$ & $l_1(v)$               & $l_2(v)$               & $l_{\infty}(v)$        & $q_2^{emp}$     \\
  \midrule
  $10 \times \hspace{2mm}5$ &    3600           & $3.608 \times 10^{-1}$ & $3.665 \times 10^{-1}$ & $5.166 \times 10^{-1}$ & -\\      
  $20 \times 10           $ &    1800           & $2.164 \times 10^{-1}$ & $2.198 \times 10^{-1}$ & $3.041 \times 10^{-1}$ & 0.7\\
  $40 \times 20           $ &    900            & $1.185 \times 10^{-1}$ & $1.203 \times 10^{-1}$ & $1.671 \times 10^{-2}$ & 0.9\\
  $80 \times 40           $ &    450            & $6.195 \times 10^{-2}$ & $6.291 \times 10^{-2}$ & $8.809 \times 10^{-2}$ & 0.9 \\  
  \bottomrule
\end{tabular}
\end{center}
\caption{Relative errors on $v $ at different resolutions,   L\"auter test case with off centered Crank Nicolson, $\theta=0.6$.}
\label{tlauterconvrate_v_tab_thet06}
\end{table}
 
\clearpage 

\subsection{Zonal flow over an isolated mountain}
\label{test5}
We have then performed numerical simulations reproducing the
test case 5 of \cite{williamson:1992}, given by a zonal flow impinging 
on an isolated mountain of conical shape.
The geostrophic balance here is broken by orographic forcing, 
which  results in the development of a planetary wave propagating all around the globe. 

%
%
%
%
%
%

\begin{figure}[htbc] 
 \centering
 \includegraphics[width=\textwidth]{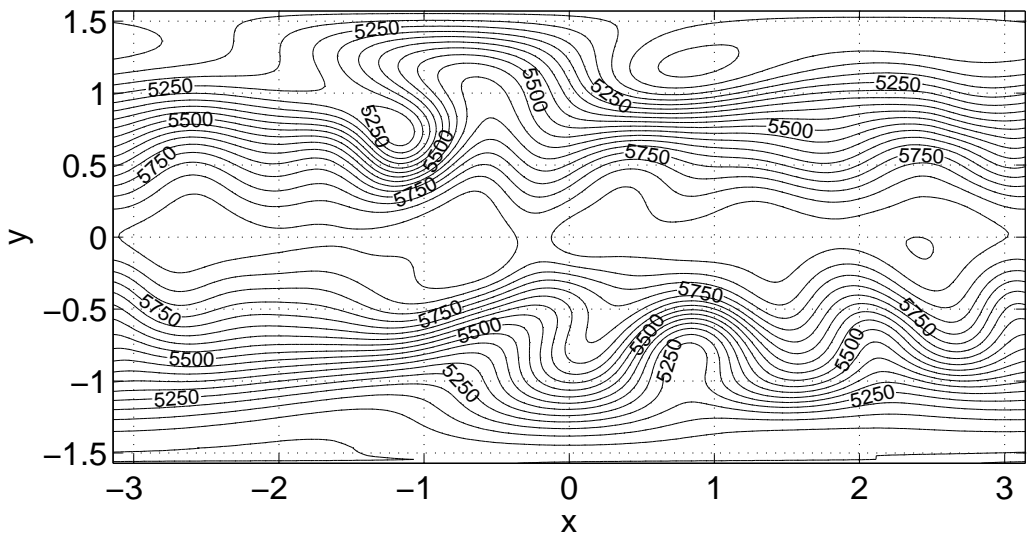}
 \caption{$h$ field after 15 days, isolated mountain wave test case, $C_{cel} \approx 58.$ Contour lines spacing is 50 m.}
\label{fig:t5_h}
\end{figure}

\begin{figure}[htbc] 
 \centering
 \includegraphics[width=\textwidth]{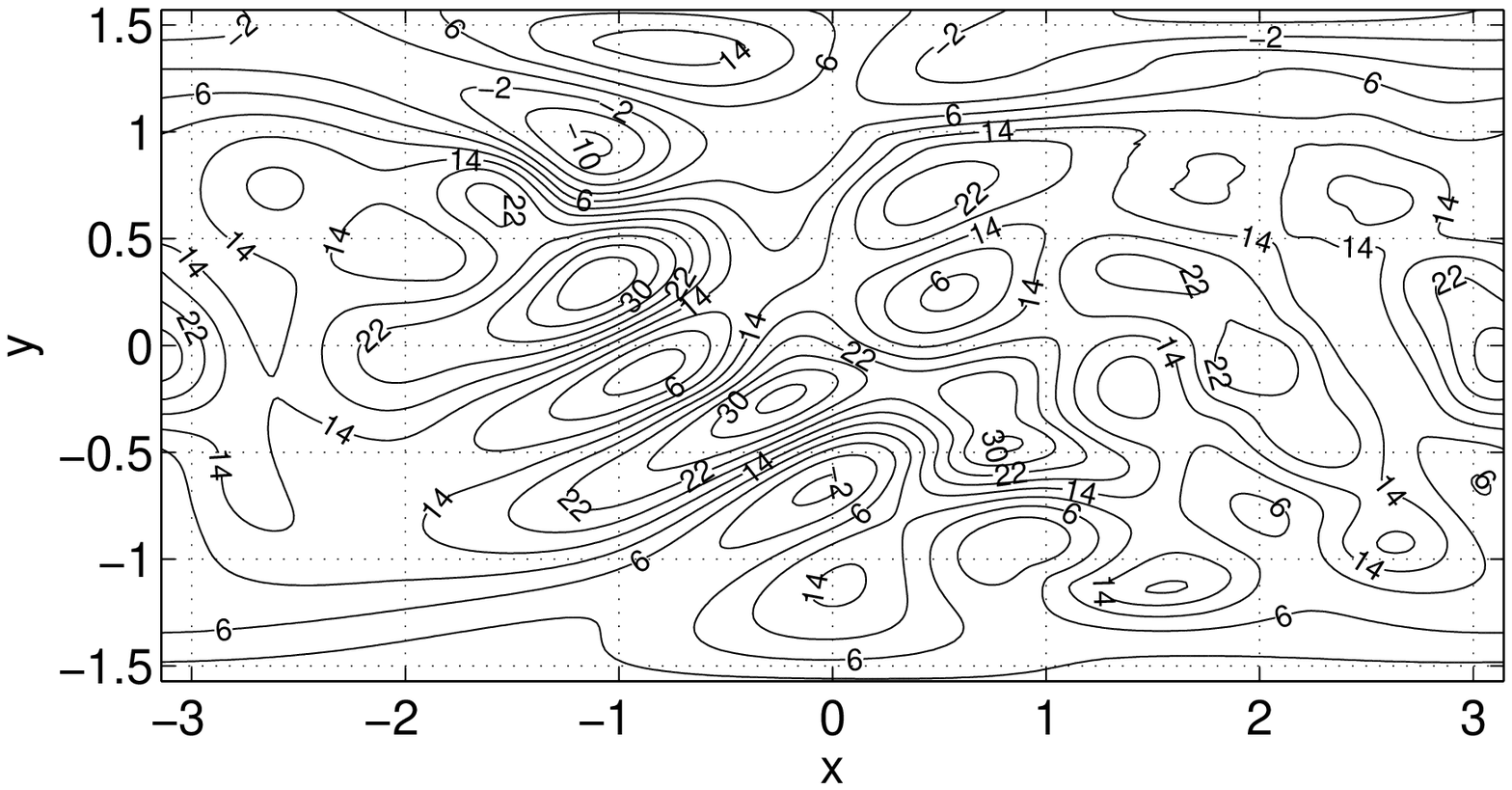}
 \caption{$u$ field after 15 days, isolated mountain wave test case, $C_{cel} \approx 58.$ Contour lines spacing is 4 m s$^{-1}$.}
\label{fig:t5_u}
\end{figure}

\begin{figure}[htbc] 
 \centering
 \includegraphics[width=\textwidth]{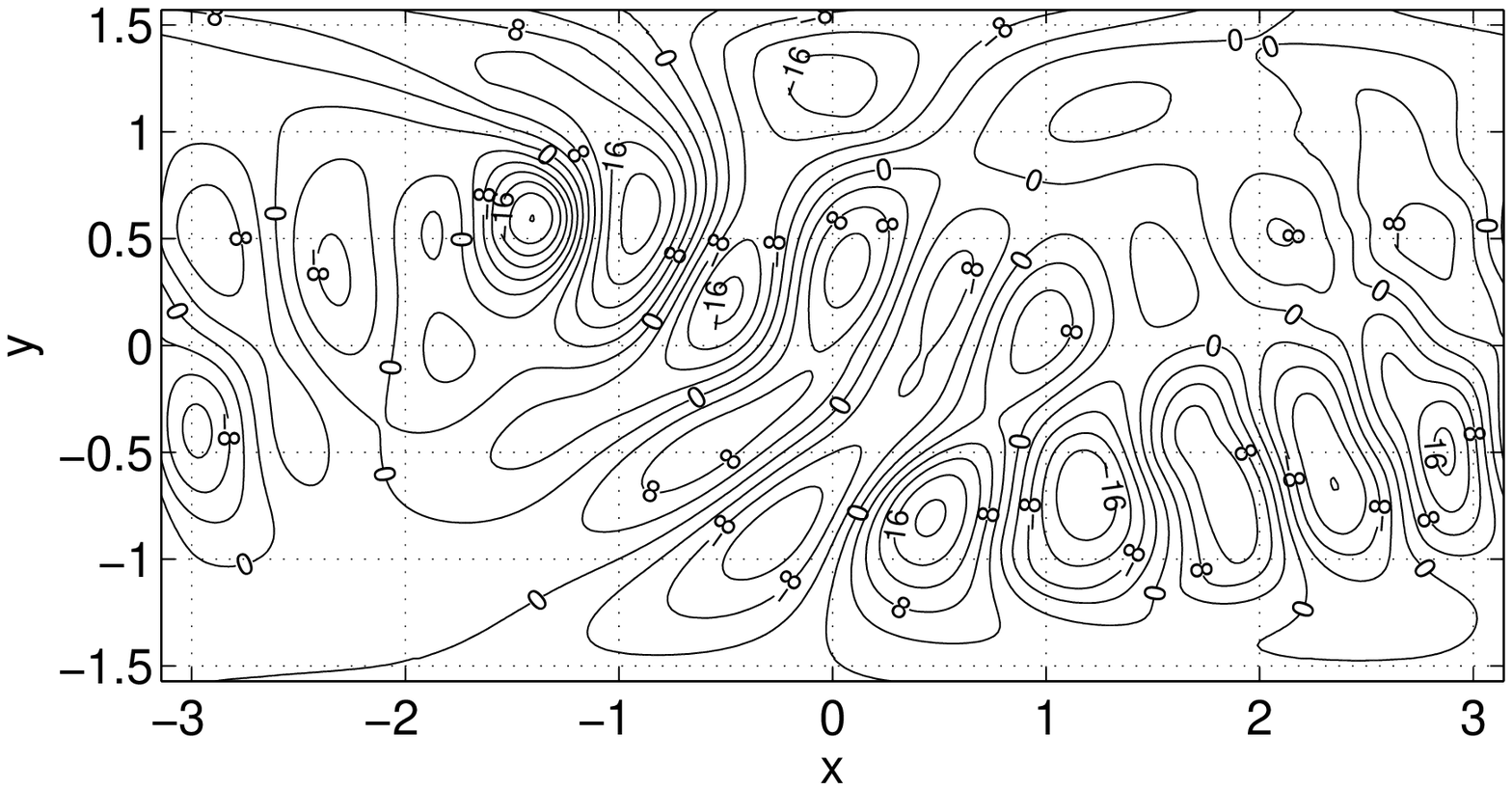}
 \caption{$v$ field after 15 days, isolated mountain wave test case, $C_{cel} \approx 58.$ Contour lines spacing is 4 m s$^{-1}$.}
\label{fig:t5_v}
\end{figure}

Plots of the fluid depth $h$ as well as of the velocity components $u$ and $v$ at 15 days 
are shown in figures \ref{fig:t5_h}-\ref{fig:t5_v}.
The resolution used corresponds to a mesh of $60 \times 30$ elements with 
$p^h = 4,$  $p^u = 5,$ and $\Delta t = 900 \ \text{s},$ giving  
a Courant number $C_{cel} \approx 58$ 
in elements close to the poles.
It can be observed that all the main features of the flow are correctly reproduced.
In particular, no significant Gibbs phenomena are detected in the vicinity of the mountain,
 even in the initial stages of the simulation.

\begin{figure}[htbc]
\subfigure[ \label{fig:test5_invariant_mass}]
  {\includegraphics[width=.55\linewidth]{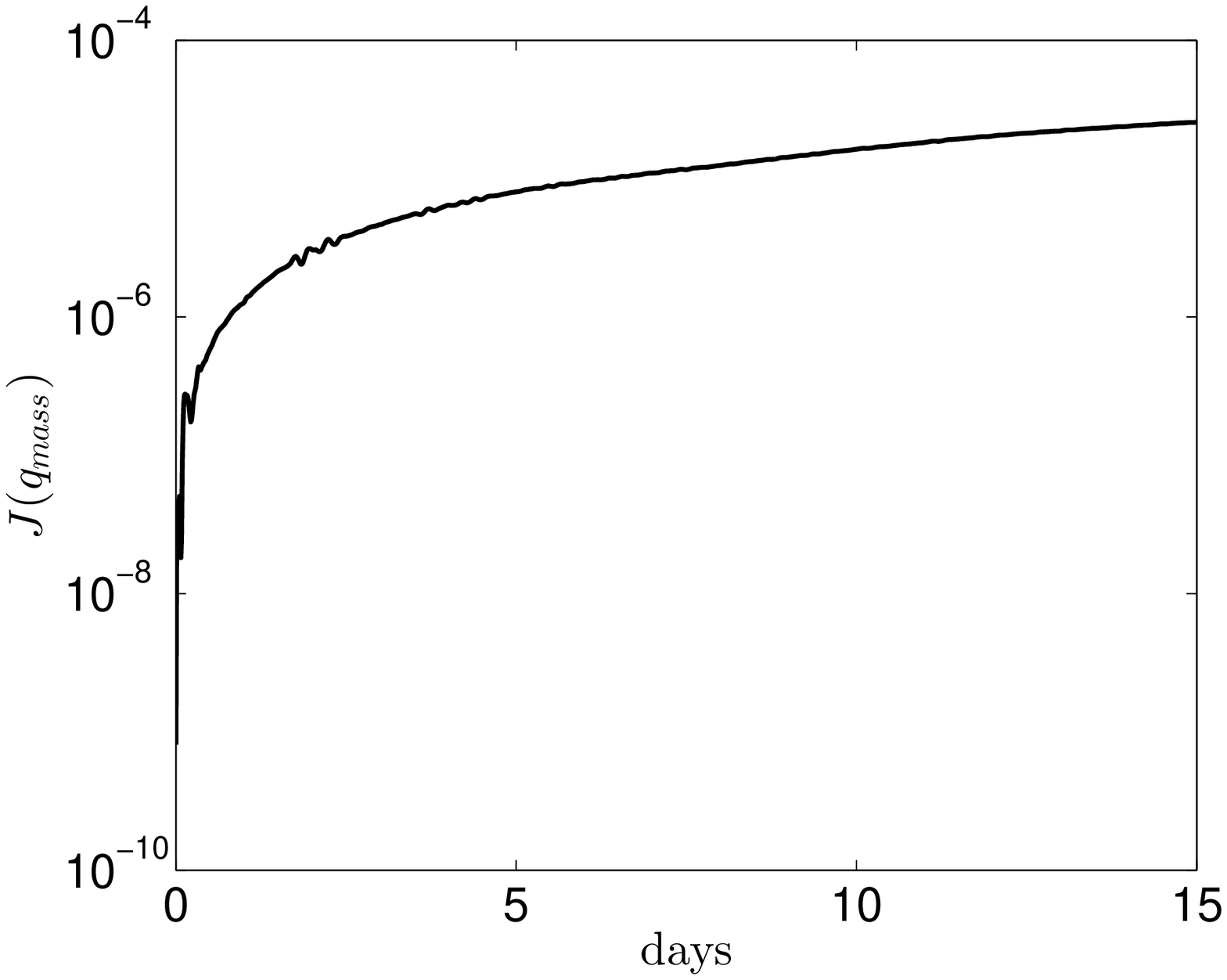}}\hfill
\subfigure[ \label{fig:test5_invariant_energy}]
  {\includegraphics[width=.55\linewidth]{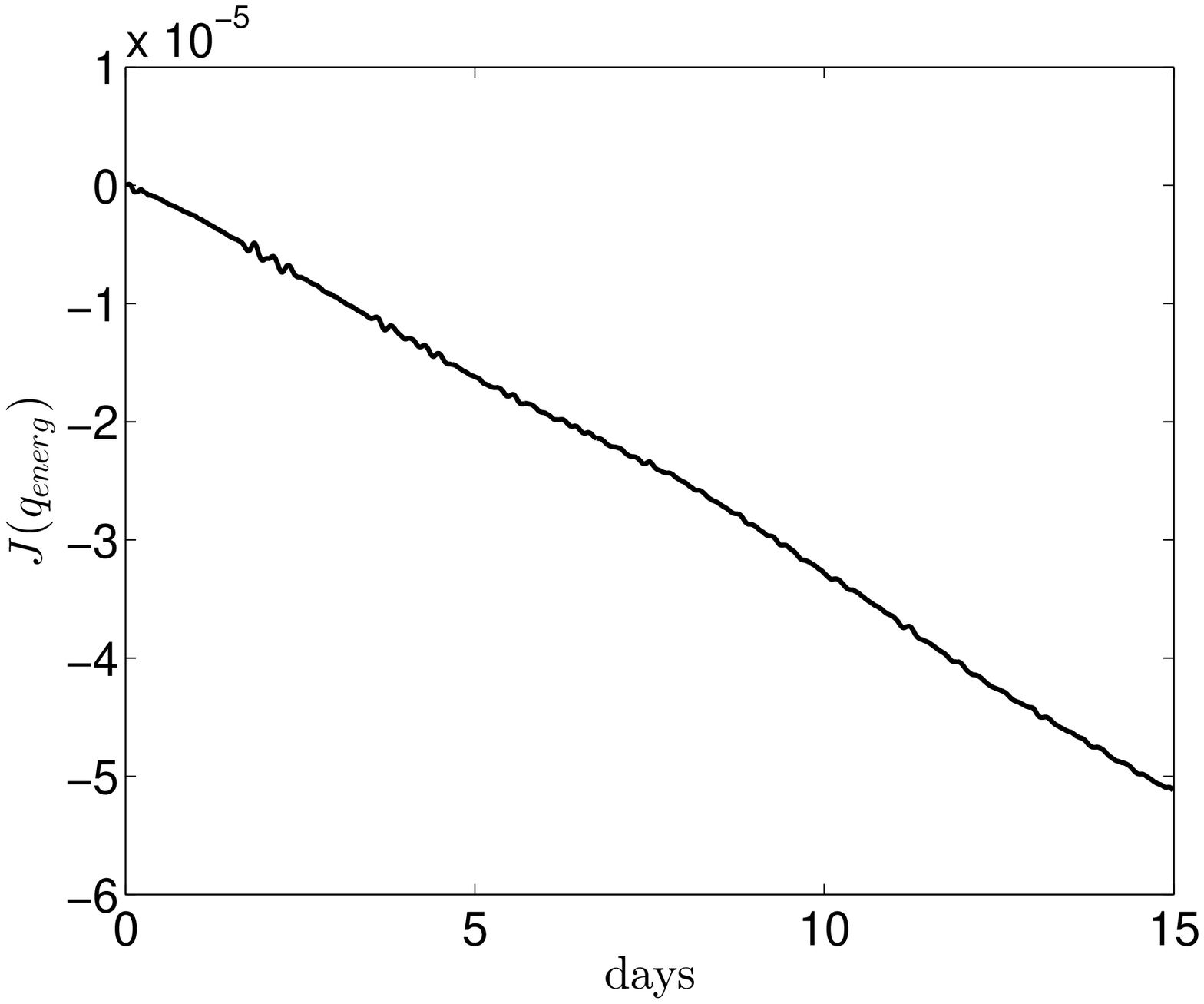}}\hfill
  \subfigure[ \label{fig:test5_invariant_enstrophy}]
  {\includegraphics[width=.55\linewidth]{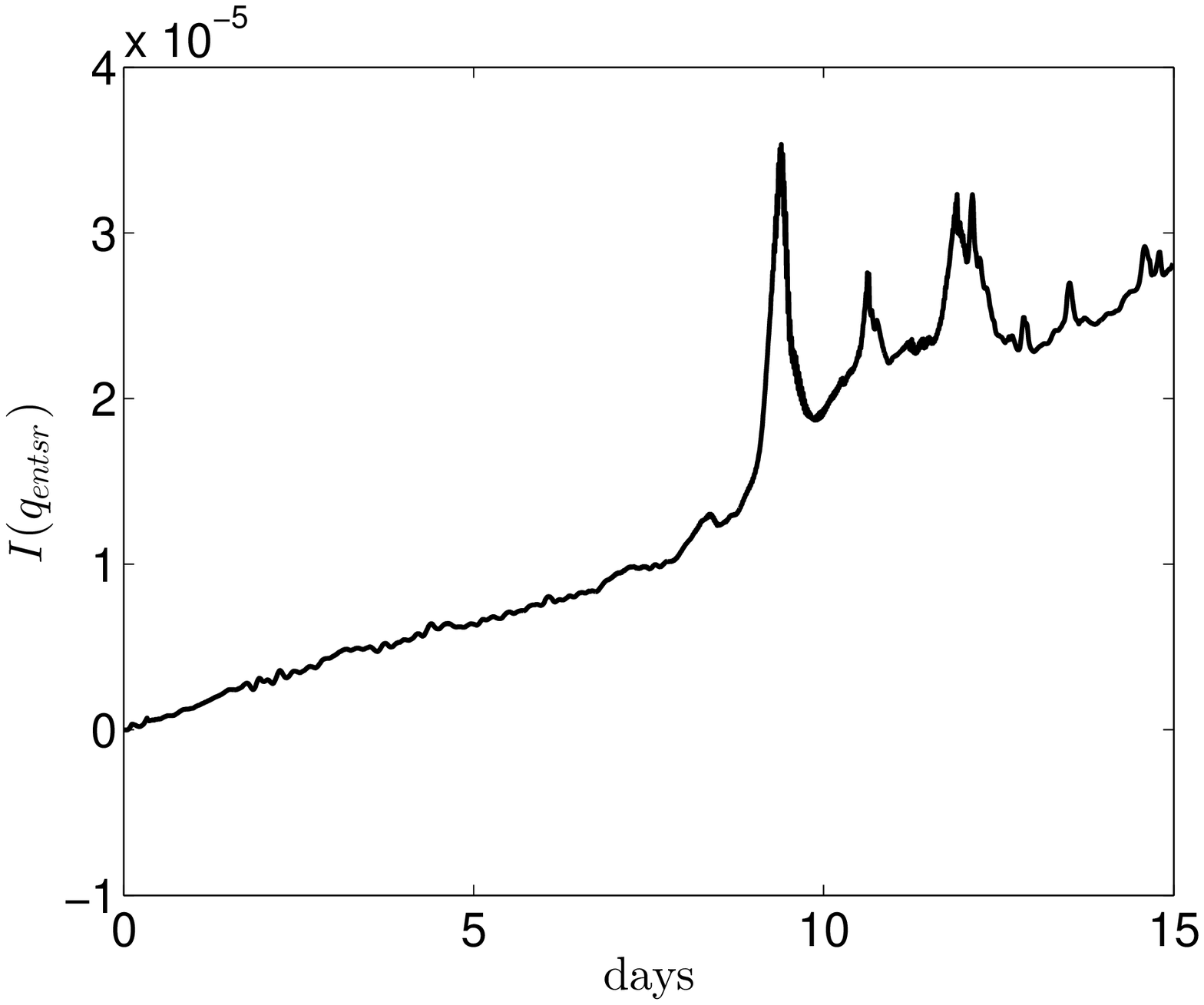}}\hfill
\caption{Integral invariants evolution, mass (a), energy (b), potential enstrophy (c),  isolated mountain wave test case, $C_{cel} \approx 58.$}
\end{figure} 



 The evolution in time of global invariants during this simulation is shown in figures \ref{fig:test5_invariant_mass},
\ref{fig:test5_invariant_energy}, \ref{fig:test5_invariant_enstrophy}, respectively.
Error norms for $h$  and $u$ at different resolutions corresponding to a $C_{cel} \approx 6$ and $p^h=p^u=3,$ have been computed at $t_f=5$ days
and are displayed in tables \ref{tab:t5convrate_h} - \ref{tab:t5convrate_u}, 
with respect to a reference solution given by the National Center for Atmospheric Research (NCAR) 
spectral model \cite{jakob:1995} at resolution T511. 
It is apparent the second order of the proposed SISLDG scheme in time.
Since, as observed in \cite{jakob:1995}, the National Center for Atmospheric Research (NCAR) 
spectral model incorporates diffusion terms in the governing equations,
while the proposed SISLDG scheme does not employ any diffusion terms nor filtering, nor smoothing of the topography,
for this test it seemed more appropriate to compute relative errors with respect to  
NCAR spectral model \cite{jakob:1995} solution
at an earlier time, $t_f=5$ days, when it can be assumed that the effects of diffusion have less impact.
Error norms for $h$  and $u$ have been computed at $t_f=15$ days at different resolutions (corresponding to a $C_{cel} \approx 7$), $p^h=p^u=3,$ 
and displayed in tables \ref{tab:t5convrate_h_15dd} - \ref{tab:t5convrate_u_15dd}.

\begin{table}[htbc] 
\begin{center}
\begin{tabular}{cccccc}
  \toprule
  $N_x \times N_y$          &  $\Delta t \text{[min]}$ & $l_1(h)$              &  $l_2(h)$              & $l_{\infty}(h)$       & $q_2^{emp}$        \\
  \midrule
  $12 \times \hspace{2mm}6$ &    20                    & $8.19 \times 10^{-4}$ &  $1.08 \times 10^{-3}$ & $5.90 \times 10^{-3}$ &   -                \\
  $24 \times 12           $ &    10                    & $1.49 \times 10^{-4}$ &  $2.08 \times 10^{-4}$ & $1.92 \times 10^{-3}$ &   2.4              \\
  $48 \times 24           $ &     5                    & $2.88 \times 10^{-5}$ &  $4.25 \times 10^{-5}$ & $8.40 \times 10^{-4}$ &   2.3              \\
  \bottomrule
\end{tabular}
\end{center}
\caption{Relative errors on $h $ at different resolutions, isolated mountain wave test case, $t_f=5$ days.}
\label{tab:t5convrate_h}
\end{table}

\begin{table}[htbc] 
\begin{center}
\begin{tabular}{cccccc}
  \toprule
  $N_x \times N_y$          & $\Delta t \text{[min]}$ & $l_1(u)$               & $l_2(u)$               & $l_{\infty}(u)$         & $q_2^{emp}$        \\
  \midrule
  $12 \times \hspace{2mm}6$ &    20                   & $4.33 \times 10^{-2}$  & $5.81 \times 10^{-2}$  & $1.39 \times 10^{-1}$   &     -              \\      
  $24 \times 12           $ &    10                   & $5.70 \times 10^{-3}$  & $7.33 \times 10^{-3}$  & $1.06 \times 10^{-1}$   &   2.9              \\
  $48 \times 24           $ &     5                   & $1.11 \times 10^{-3}$  & $1.72 \times 10^{-3}$  & $1.56 \times 10^{-2}$   &   2.2              \\
  \bottomrule
\end{tabular}
\end{center}
\caption{Relative errors on $u $ at different resolutions, isolated mountain wave test case, $t_f=5$ days.}
\label{tab:t5convrate_u}
\end{table}

\begin{table}[htbc] 
\begin{center}
\begin{tabular}{cccccc}
  \toprule
  $N_x \times N_y$          &  $\Delta t \text{[s]}$ & $l_1(h)$              &  $l_2(h)$              & $l_{\infty}(h)$       & $q_2^{emp}$        \\
  \midrule
  $12 \times \hspace{2mm}6$ &    1500                & $2.34 \times 10^{-3}$ &  $2.92 \times 10^{-3}$ & $1.49 \times 10^{-2}$ &   -                \\
  $24 \times 12           $ &     750                & $5.99 \times 10^{-4}$ &  $7.72 \times 10^{-4}$ & $3.87 \times 10^{-3}$ &   1.9              \\
  $48 \times 24           $ &     375                & $2.00 \times 10^{-4}$ &  $2.74 \times 10^{-4}$ & $1.87 \times 10^{-3}$ &   1.5              \\
  \bottomrule
\end{tabular}
\end{center}
\caption{Relative errors on $h $ at different resolutions, isolated mountain wave test case, $t_f=15$ days.}
\label{tab:t5convrate_h_15dd}
\end{table}

\begin{table}[htbc] 
\begin{center}
\begin{tabular}{cccccc}
  \toprule
  $N_x \times N_y$          & $\Delta t \text{[s]}$ & $l_1(u)$               & $l_2(u)$               & $l_{\infty}(u)$         & $q_2^{emp}$        \\
  \midrule
  $12 \times \hspace{2mm}6$ &    1500               & $1.12 \times 10^{-1}$  & $1.29 \times 10^{-1}$  & $2.97 \times 10^{-1}$   &     -              \\      
  $24 \times 12           $ &     750               & $2.09 \times 10^{-2}$  & $2.37 \times 10^{-2}$  & $5.73 \times 10^{-2}$   &   2.4              \\
  $48 \times 24           $ &     375               & $6.37 \times 10^{-3}$  & $7.92 \times 10^{-3}$  & $3.39 \times 10^{-2}$   &   1.6              \\
  \bottomrule
\end{tabular}
\end{center}
\caption{Relative errors on $u $ at different resolutions, isolated mountain wave test case, $t_f=15$ days.}
\label{tab:t5convrate_u_15dd}
\end{table}

\begin{figure}[htbc] 
 \centering
 \includegraphics[width= \textwidth]{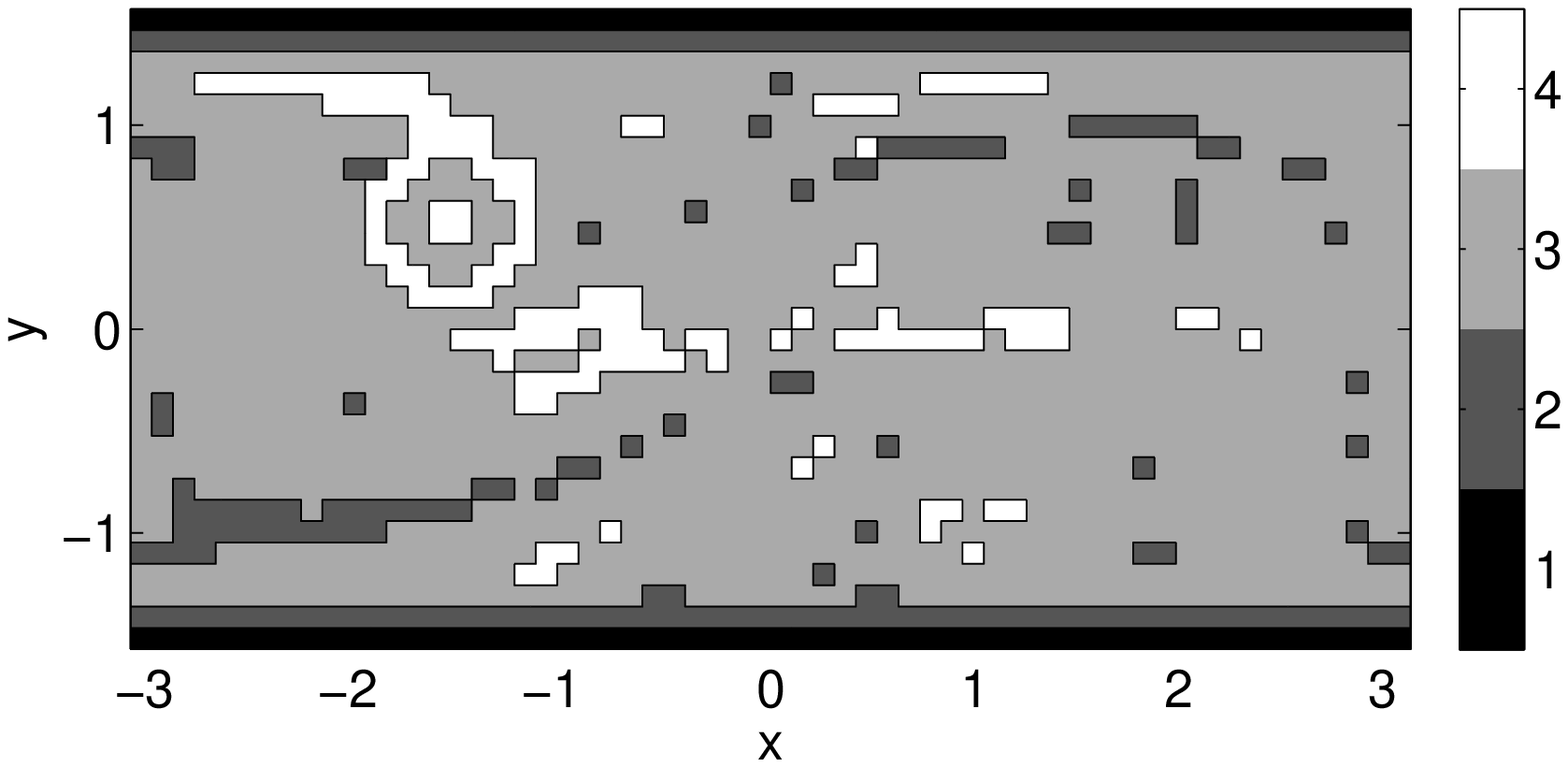}
 \caption{Statically and dynamically adapted local $p^h$ distribution at 15 days, isolated mountain wave test case.}
\label{fig:t5_ph}
\end{figure}

Finally the mountain wave test case has been run on the same mesh of $60 \times 30 $ elements, $\Delta t = 900 $ s,
with either static or  static plus dynamic adaptivity. The tolerance $\epsilon$ for the dynamic adaptivity \cite{tumolo:2013}
has been set to $\epsilon = 10^{-2}$. Results are reported 
in terms of error norms with respect to a nonadaptive solution at the maximum uniform resolution and in terms of efficiency gain,
measured through the saving of number of linear solver iterations per time-step $\Delta_{iter}^{average} $ as well as
through the saving of number of degrees of freedom actually used per timestep $\Delta_{dof}^{average} $; these results
are summarized in tables \ref{t5_adaptivity_h_tab} - \ref{t5_adaptivity_v_tab}:
the use of static adaptivity only resulted in $\Delta_{iter}^{average} \approx 10.7\%$ and $\Delta_{dof}^{average} \approx 88\%,$
while the use of both static and dynamic adaptivity led to $\Delta_{iter}^{average} \approx 13\%$ and $\Delta_{dof}^{average} \approx 45\%.$
The distribution of the statically and dynamically adapted local polynomial degree used 
to represent the solution after 15 days is shown in figure \ref{fig:t5_ph}. 
It can be noticed how, even after 15 days, higher polynomial degrees are still automatically concentrated around the location of the mountain. 

\begin{table}[htbc] 
\begin{center}
\begin{tabular}{cccc}
  \toprule
  $\text{adaptivity}$       & $l_1(h)$               & $l_2(h)$               & $l_{\infty}(h)$        \\
  \midrule
  $\text{static}$           & $1.415 \times 10^{-4}$ & $3.314 \times 10^{-4}$ & $2.117 \times 10^{-3}$ \\      
  $\text{static + dynamic}$ & $1.660 \times 10^{-4}$ & $3.419 \times 10^{-4}$ & $2.038 \times 10^{-3}$ \\  
  \bottomrule
\end{tabular}
\end{center}
\caption{Relative errors between (statically and statically plus dynamically) adaptive and nonadaptive solution for isolated mountain wave test case, $h $ field.}
\label{t5_adaptivity_h_tab}
\end{table}

\begin{table}[htbc] 
\begin{center}
\begin{tabular}{cccc}
  \toprule
  $\text{adaptivity}$       & $l_1(u)$                & $l_2(u)$               & $l_{\infty}(u)$       \\
  \midrule
  $\text{static}$           & $1.289 \times 10^{-2}$ & $3.275 \times 10^{-2}$ & $1.524 \times 10^{-1}$ \\      
  $\text{static + dynamic}$ & $1.509 \times 10^{-2}$ & $3.309 \times 10^{-2}$ & $1.475 \times 10^{-1}$ \\  
  \bottomrule
\end{tabular}
\end{center}
\caption{Relative errors between (statically and statically plus dynamically) adaptive and nonadaptive solution for isolated mountain wave test case,  $u $ field.}
\label{t5_adaptivity_u_tab}
\end{table}

\begin{table}[htbc] 
\begin{center}
\begin{tabular}{cccc}
  \toprule
  $\text{adaptivity}$       & $l_1(v)$                & $l_2(v)$               & $l_{\infty}(v)$       \\
  \midrule
  $\text{static}$           & $2.501 \times 10^{-2}$ & $6.824 \times 10^{-2}$ & $6.597 \times 10^{-1}$ \\      
  $\text{static + dynamic}$ & $2.833 \times 10^{-2}$ & $7.019 \times 10^{-2}$ & $6.975 \times 10^{-1}$ \\  
  \bottomrule
\end{tabular}
\end{center}
\caption{Relative errors between (statically and statically plus dynamically) adaptive and nonadaptive solution for isolated mountain wave test case,  $v $ field.}
\label{t5_adaptivity_v_tab}
\end{table}

\clearpage

\subsection{Rossby-Haurwitz wave}
\label{rossby}

\begin{figure}[htbc] 
 \centering
 \includegraphics[width=\textwidth]{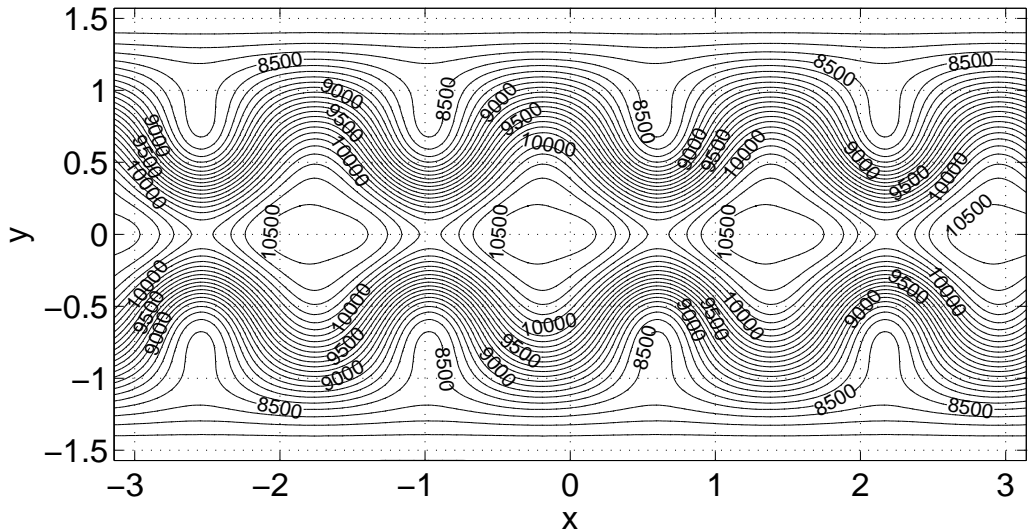}
 \caption{ $h$ field after 15 days, Rossby-Haurwitz wave test case, $C_{cel} \approx 83.$ Contour lines spacing is 100 m.}
\label{fig:t6_h}
\end{figure}

\begin{figure}[htbc] 
 \centering
 \includegraphics[width=\textwidth]{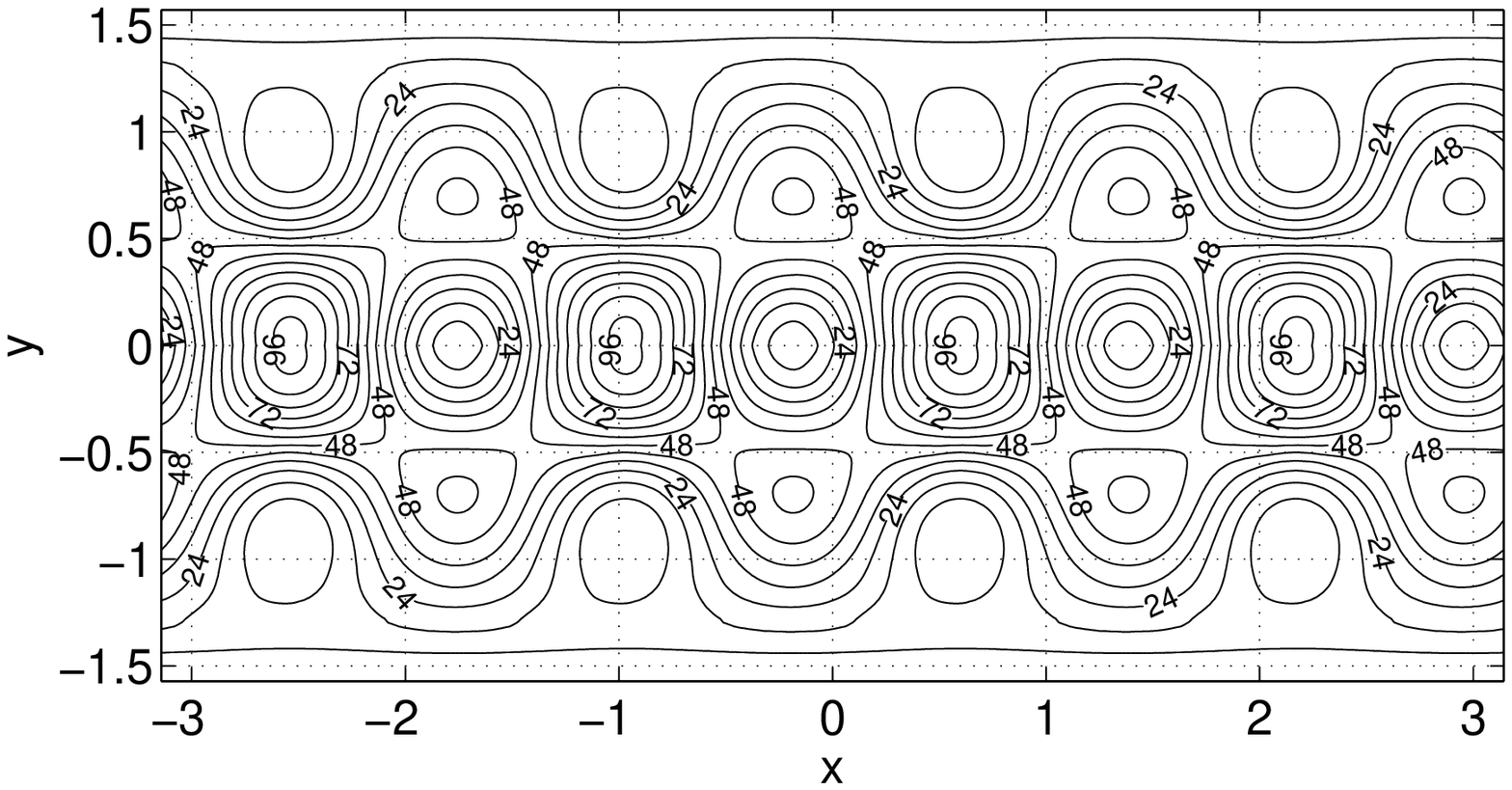}
 \caption{$u$ field after 15 days, Rossby-Haurwitz wave test case, $C_{cel} \approx 83.$ Contour lines spacing is 8 m s$^{-1}$.}
\label{fig:t6_u}
\end{figure}

\begin{figure}[htbc] 
 \centering
 \includegraphics[width=\textwidth]{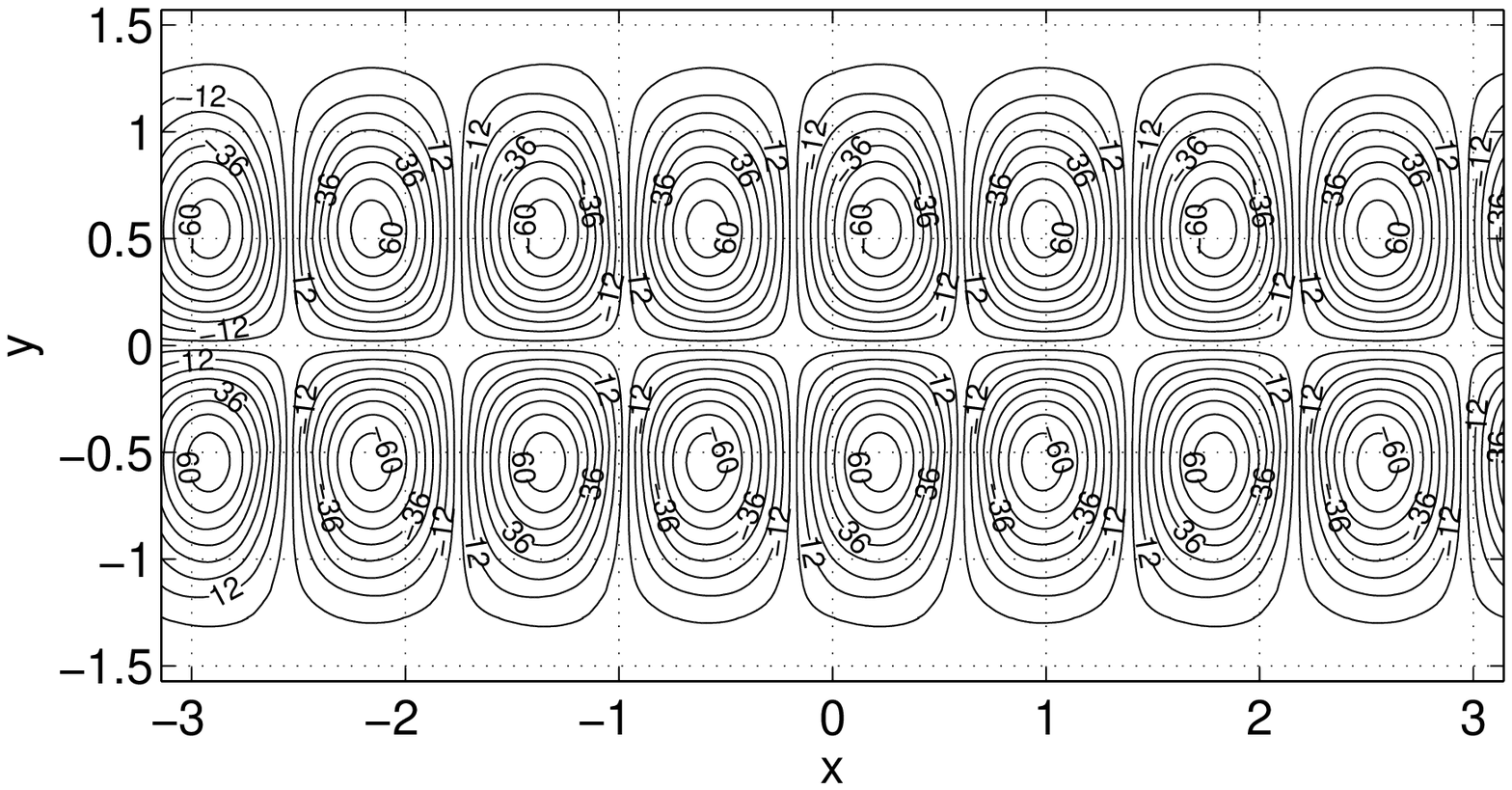}
 \caption{$v$ field after 15 days, Rossby-Haurwitz wave test case, $C_{cel} \approx 83.$ Contour lines spacing is 8 m s$^{-1}$.}
\label{fig:t6_v}
\end{figure}

We have then considered  test case 6 of \cite{williamson:1992}, where the initial datum consists
of a Rossby-Haurwitz wave of wave number 4. This case  actually concerns a   solution of
the nondivergent barotropic vorticity equation, that is  not an exact solution of the system (\ref{continuityeq}) - (\ref{vectmomentumeq}).
For a discussion about the stability of this profile as a solution of (\ref{continuityeq}) - (\ref{vectmomentumeq}) see \cite{thuburn:2000}.
Plots of the fluid depth $h$ as well as of the velocity components $u$ and $v$ at 15 days 
are shown in figures \ref{fig:t6_h}-\ref{fig:t6_v}.
The resolution used corresponds to a mesh of $64 \times 32$ 
elements with $p^h = 4,$  $p^u = 5,$ and $\Delta t = 900 \ \text{s},$ giving  
a Courant number $C_{cel} \approx 83$ 
in elements close to poles.
It can be observed that all the main features of the flow are correctly reproduced.

\begin{figure}[htbc]
\subfigure[ \label{fig:test6_invariant_mass}]
  {\includegraphics[width=.55\linewidth]{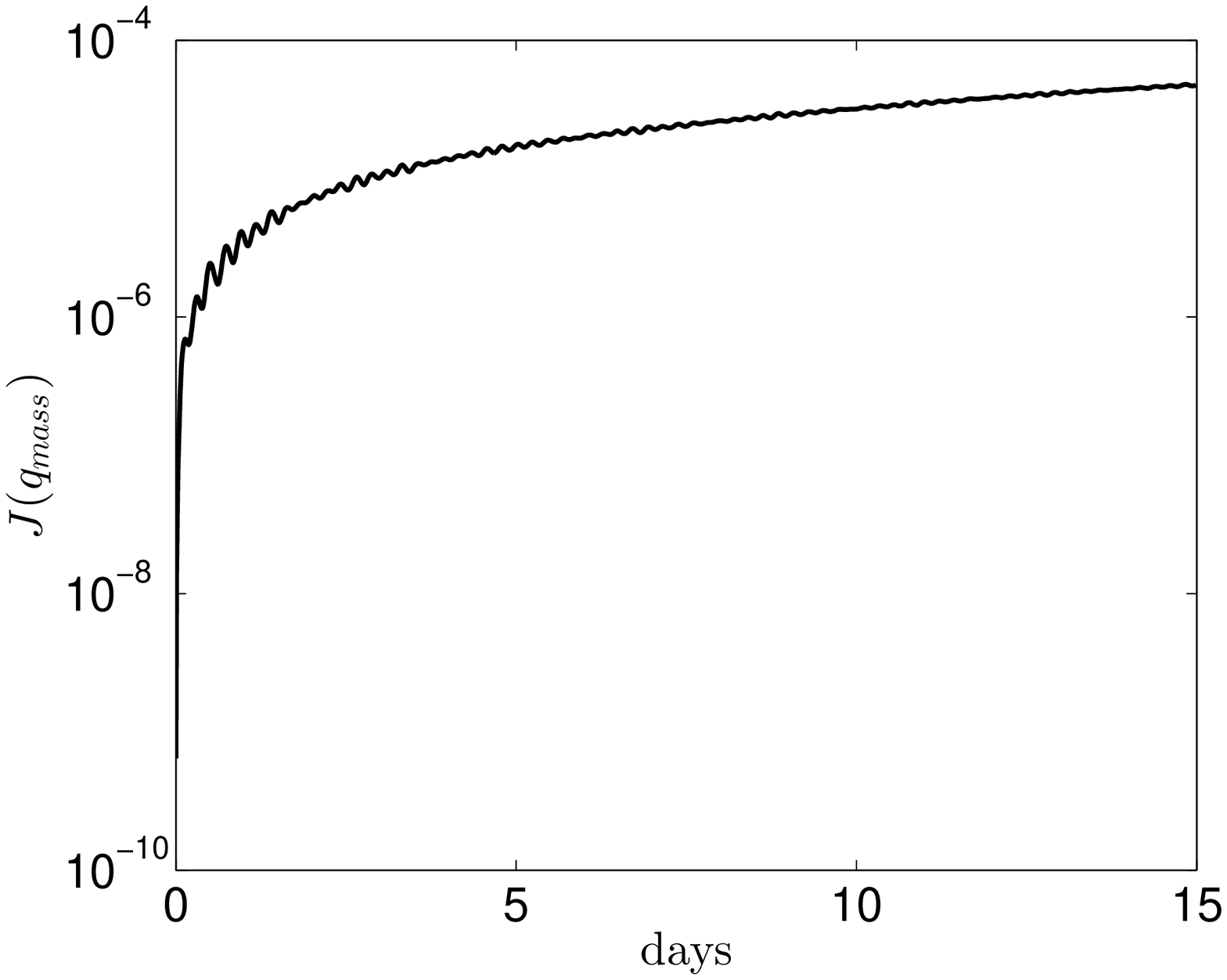}}\hfill
\subfigure[ \label{fig:test6_invariant_energy}]
  {\includegraphics[width=.55\linewidth]{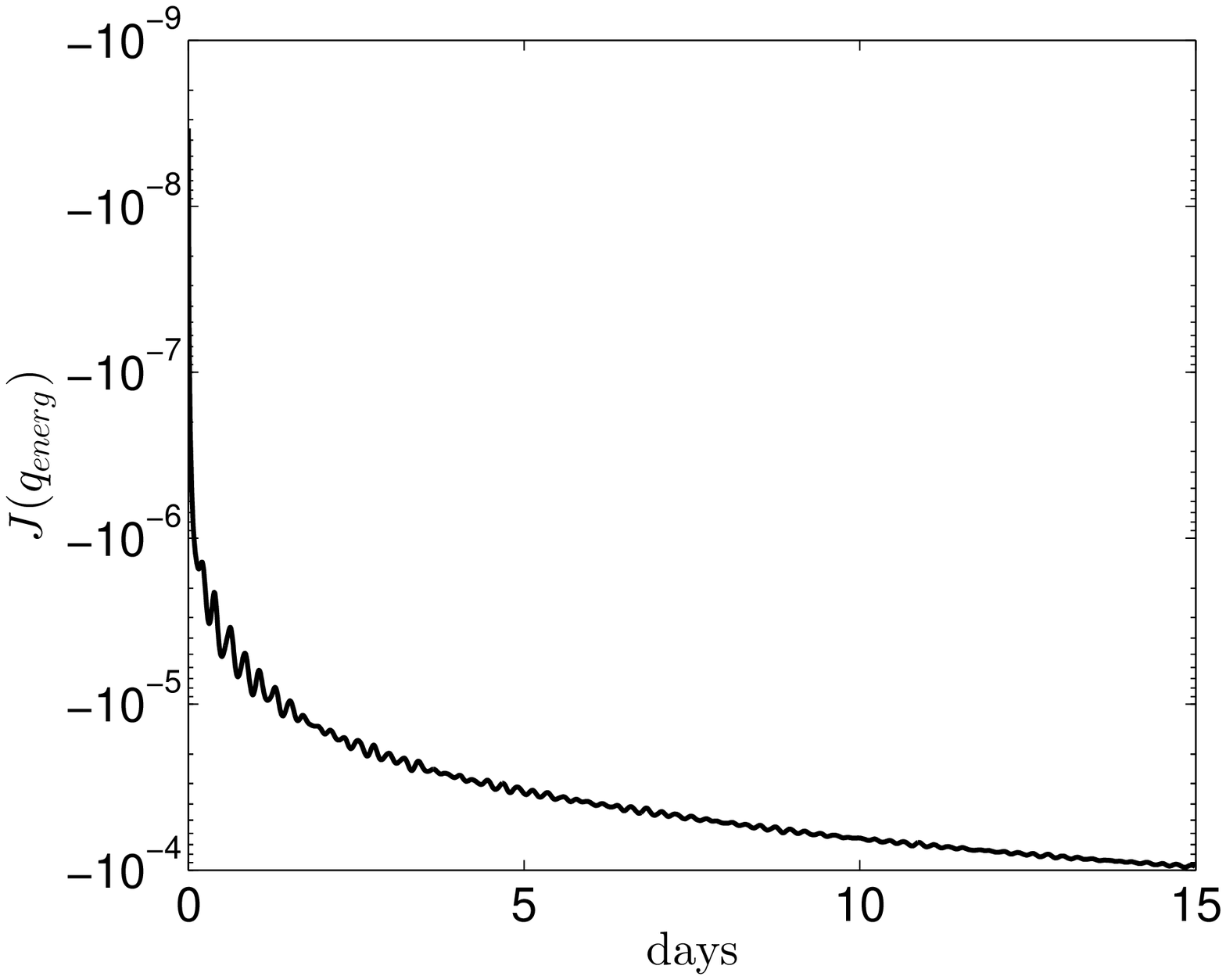}}\hfill
  \subfigure[ \label{fig:test6_invariant_enstrophy}]
  {\includegraphics[width=.55\linewidth]{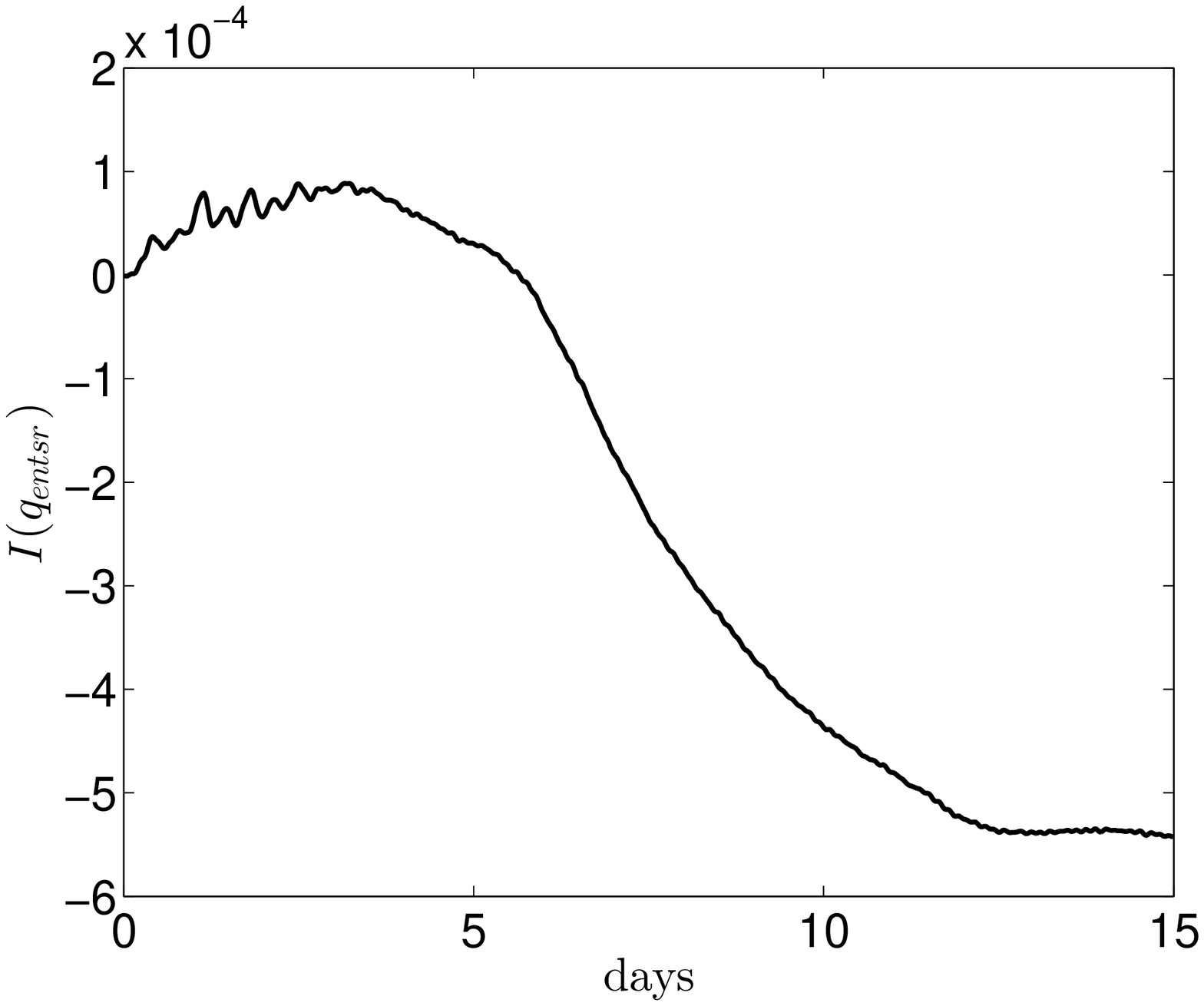}}\hfill
\caption{Integral invariants evolution, mass (a), energy (b), potential enstrophy (c), Rossby-Haurwitz wave test case, $C_{cel} \approx 83.$}
\end{figure}

The evolution in time of global invariants during this simulation is shown in figures \ref{fig:test6_invariant_mass},
\ref{fig:test6_invariant_energy}, \ref{fig:test6_invariant_enstrophy}, respectively.
Error norms for $h$  and $u$ at different resolutions, corresponding to a $C_{cel} \approx 32 $ and $ p^h=4, p^u=5 $ have been computed at $t_f=15$ days
and are displayed in tables \ref{tab:t6convrate_h} - \ref{tab:t6convrate_u}, 
with respect to a reference solution given by the National Center for Atmospheric Research (NCAR) 
spectral model \cite{jakob:1995} at resolution T511. 
It is apparent the second order of the proposed SISLDG scheme in time.
Unlike the NCAR spectral model, the proposed SISLDG scheme does not employ any explicit numerical diffusion.

\begin{table}[htbc] 
\begin{center}
\begin{tabular}{cccccc}
  \toprule
  $N_x \times N_y$          &  $\Delta t \text{[min]}$ & $l_1(h)$              &  $l_2(h)$              & $l_{\infty}(h)$       & $q_2^{emp}$        \\
  \midrule
  $10 \times \hspace{2mm}5$ &    60                    & $2.92 \times 10^{-2}$ &  $3.82 \times 10^{-2}$ & $6.75 \times 10^{-2}$ &   -                \\
  $20 \times 10           $ &    30                    & $5.50 \times 10^{-3}$ &  $6.80 \times 10^{-3}$ & $1.11 \times 10^{-2}$ &   2.4              \\
  $40 \times 20           $ &    15                    & $1.40 \times 10^{-3}$ &  $1.80 \times 10^{-3}$ & $3.20 \times 10^{-3}$ &   2.0              \\
  \bottomrule
\end{tabular}
\end{center}
\caption{Relative errors on $h $ at different resolutions, Rossby-Haurwitz wave test case.}
\label{tab:t6convrate_h}
\end{table}

\begin{table}[htbc] 
\begin{center}
\begin{tabular}{cccccc}
  \toprule
  $N_x \times N_y$          & $\Delta t \text{[min]}$ & $l_1(u)$               & $l_2(u)$               & $l_{\infty}(u)$         & $q_2^{emp}$        \\
  \midrule
  $10 \times \hspace{2mm}5$ &    60                   & $4.065 \times 10^{-1}$ & $3.775 \times 10^{-1}$ & $2.305 \times 10^{-1}$  &     -              \\      
  $20 \times 10           $ &    30                   & $7.79 \times 10^{-2}$  & $7.33 \times 10^{-2}$  & $5.67 \times 10^{-2}$   &   2.4              \\
  $40 \times 20           $ &    15                   & $2.04 \times 10^{-2}$  & $1.95 \times 10^{-2}$  & $1.76 \times 10^{-2}$   &   1.9              \\
  \bottomrule
\end{tabular}
\end{center}
\caption{Relative errors on $u $ at different resolutions, Rossby-Haurwitz wave test case.}
\label{tab:t6convrate_u}
\end{table}

\begin{figure}[htbc] 
 \centering
 \includegraphics[width= \textwidth]{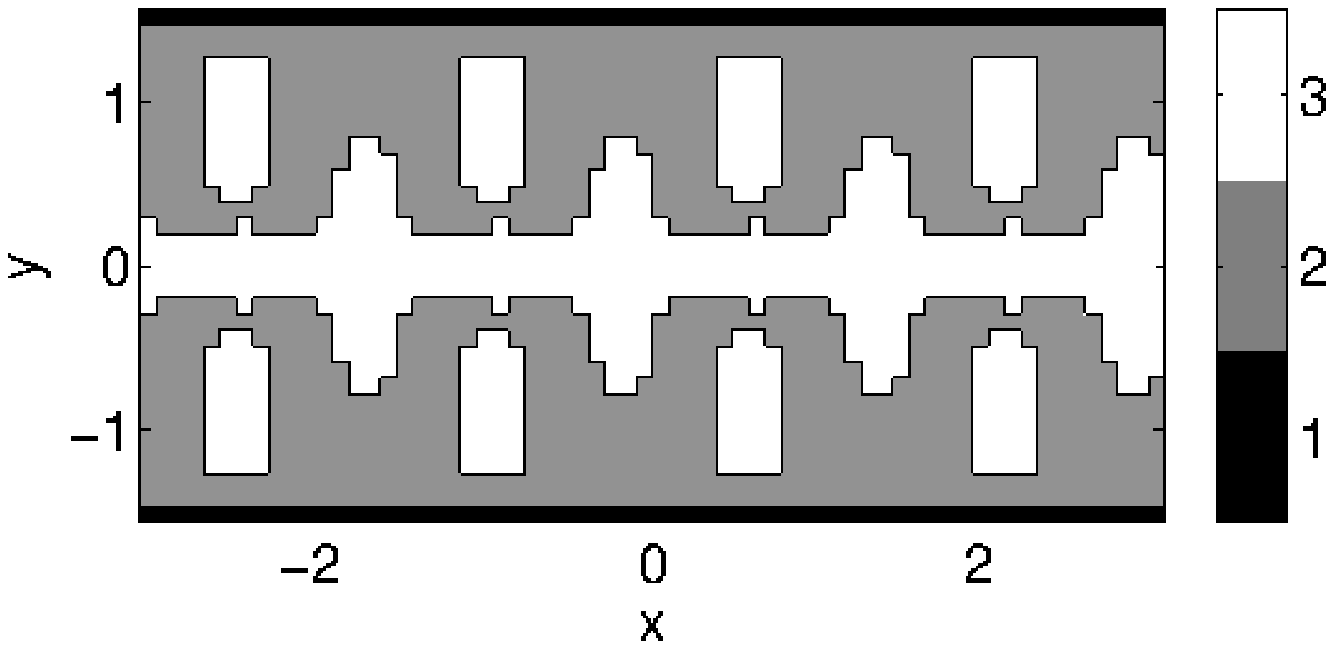}
 \caption{Statically and dynamically adapted local $p^h$ distribution at 15 days, Rossby-Haurwitz test case.}
\label{fig:t6_ph}
\end{figure}

Finally, the Rossby-Haurwitz wave test case has been run on the same mesh of $64 \times 32 $ elements, $\Delta t = 900 $ s,
with either static or  static plus dynamic adaptivity. The tolerance $\epsilon$ for the dynamic adaptivity \cite{tumolo:2013}
has been set to $\epsilon = 5 \times 10^{-2}$.  Results are reported
in terms of error norms with respect to a nonadaptive solution at the maximum uniform resolution and in terms of efficiency gain,
measured through the saving of number of linear solver iterations per time-step $\Delta_{iter}^{average} $ as well as
through the saving of number of degrees of freedom actually used per timestep $\Delta_{dof}^{average} $; these results
are summarized in tables \ref{t6_adaptivity_h_tab} - \ref{t6_adaptivity_v_tab}:
the use of static adaptivity only resulted in $\Delta_{iter}^{average} \approx 10.7\%$ and $\Delta_{dof}^{average} \approx 88\%,$
while the use of both static and dynamic adaptivity led to $\Delta_{iter}^{average} \approx 13\%$ and $\Delta_{dof}^{average} \approx 45\%.$
The distribution of the statically and dynamically adapted local polynomial degree used 
to represent the solution after 15 days is shown in figure \ref{fig:t6_ph}. 
It  can be noticed how, even after 15 days, and even if the maximum allowed $p^h$ is 4, the use of the adaptivity criterion with 
$\epsilon = 5 \times 10^{-2}$ leads to the use of at most cubic polynomials for the local representation of $h.$ 

\begin{table}[htbc] 
\begin{center}
\begin{tabular}{cccc}
  \toprule
  $\text{adaptivity}$       &    $l_1(h)$                & $l_2(h)$               & $l_{\infty}(h)$        \\
  \midrule
  $\text{static}$           &    $2.182 \times 10^{-4}$  & $3.434 \times 10^{-4}$ & $2.856 \times 10^{-4}$ \\      
  $\text{static + dynamic}$ &    $2.358 \times 10^{-3}$  & $2.963 \times 10^{-3}$ & $5.157 \times 10^{-3}$ \\  
  \bottomrule
\end{tabular}
\end{center}
\caption{Relative errors between (statically and statically plus dynamically) adaptive and nonadaptive solution for Rossby-Haurwitz wave test case, $h $ field.}
\label{t6_adaptivity_h_tab}
\end{table}

\begin{table}[htbc] 
\begin{center}
\begin{tabular}{cccc}
  \toprule
  $\text{adaptivity}$       &  $l_1(u)$               & $l_2(u)$               & $l_{\infty}(u)$        \\
  \midrule
  $\text{static}$           &  $7.041 \times 10^{-3}$ & $1.236 \times 10^{-2}$ & $2.834 \times 10^{-2}$ \\      
  $\text{static + dynamic}$ &  $3.639 \times 10^{-2}$ & $3.387 \times 10^{-2}$ & $2.678 \times 10^{-2}$ \\  
  \bottomrule
\end{tabular}
\caption{Relative errors between (statically and statically plus dynamically) adaptive and nonadaptive solution for Rossby-Haurwitz wave test case, $u $ field.}
\label{t6_adaptivity_u_tab}
\end{center}
\end{table}

\begin{table}[htbc] 
\begin{center}
\begin{tabular}{cccc}
  \toprule
  $\text{adaptivity}$       & $l_1(v)$               & $l_2(v)$               & $l_{\infty}(v)$        \\
  \midrule
  $\text{static}$           & $3.158 \times 10^{-3}$ & $3.250 \times 10^{-3}$ & $1.148 \times 10^{-2}$ \\      
  $\text{static + dynamic}$ & $2.723 \times 10^{-2}$ & $2.432 \times 10^{-2}$ & $2.646 \times 10^{-2}$ \\  
  \bottomrule
\end{tabular}
\caption{Relative errors between (statically and statically plus dynamically) adaptive and nonadaptive solution for Rossby-Haurwitz wave test case, $v $ field.}
\label{t6_adaptivity_v_tab}
\end{center}
\end{table}

\clearpage

\subsection{Nonhydrostatic inertia gravity waves}
\label{nh_igw}
\begin{figure}[htbc] 
 \centering
 \includegraphics[width= \textwidth]{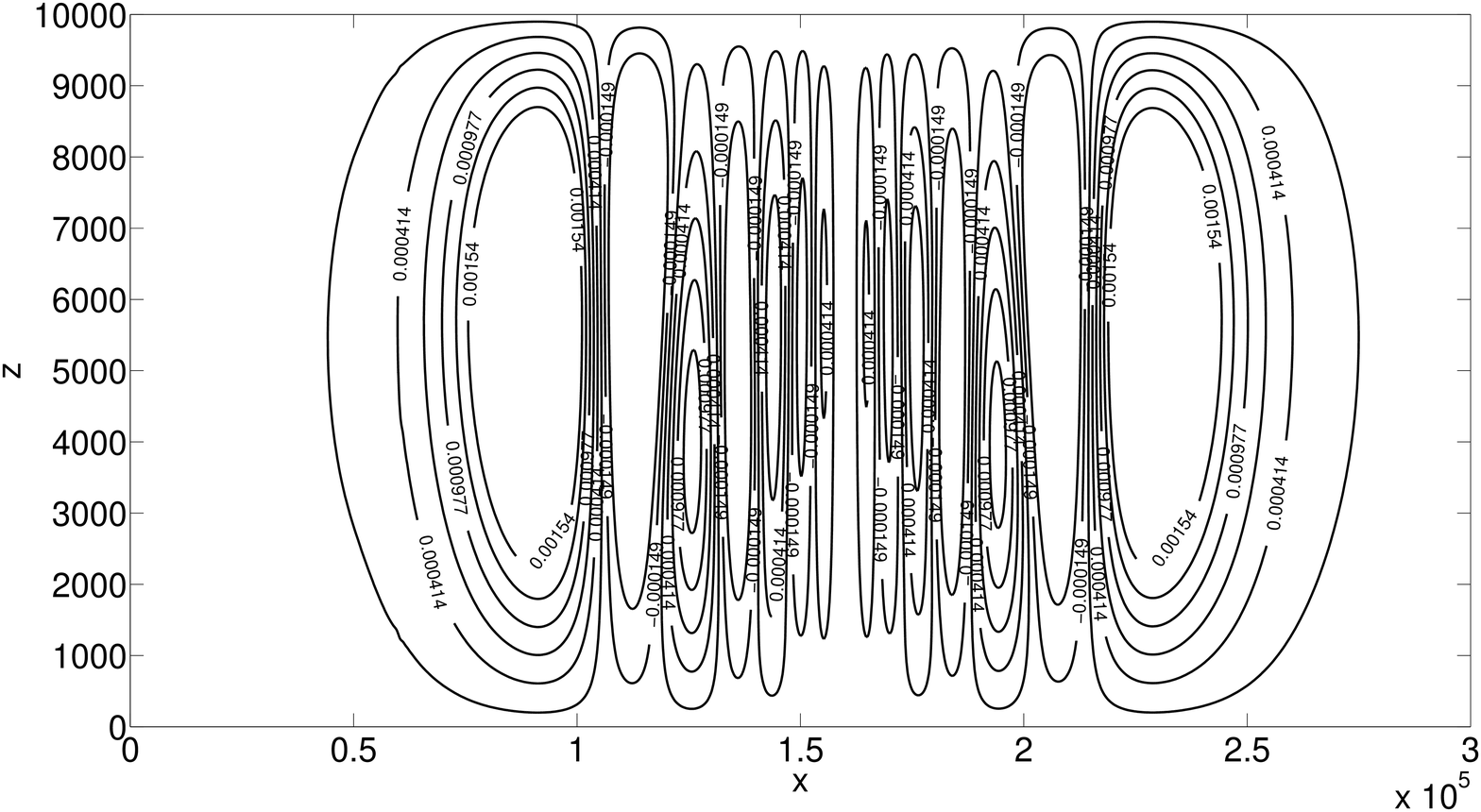}
 \caption{Contours of perturbation potential temperature in the internal gravity wave test.}
\label{fig:igw_pottemp}
\end{figure}

\begin{figure}[htbc] 
 \centering
 \includegraphics[width= \textwidth]{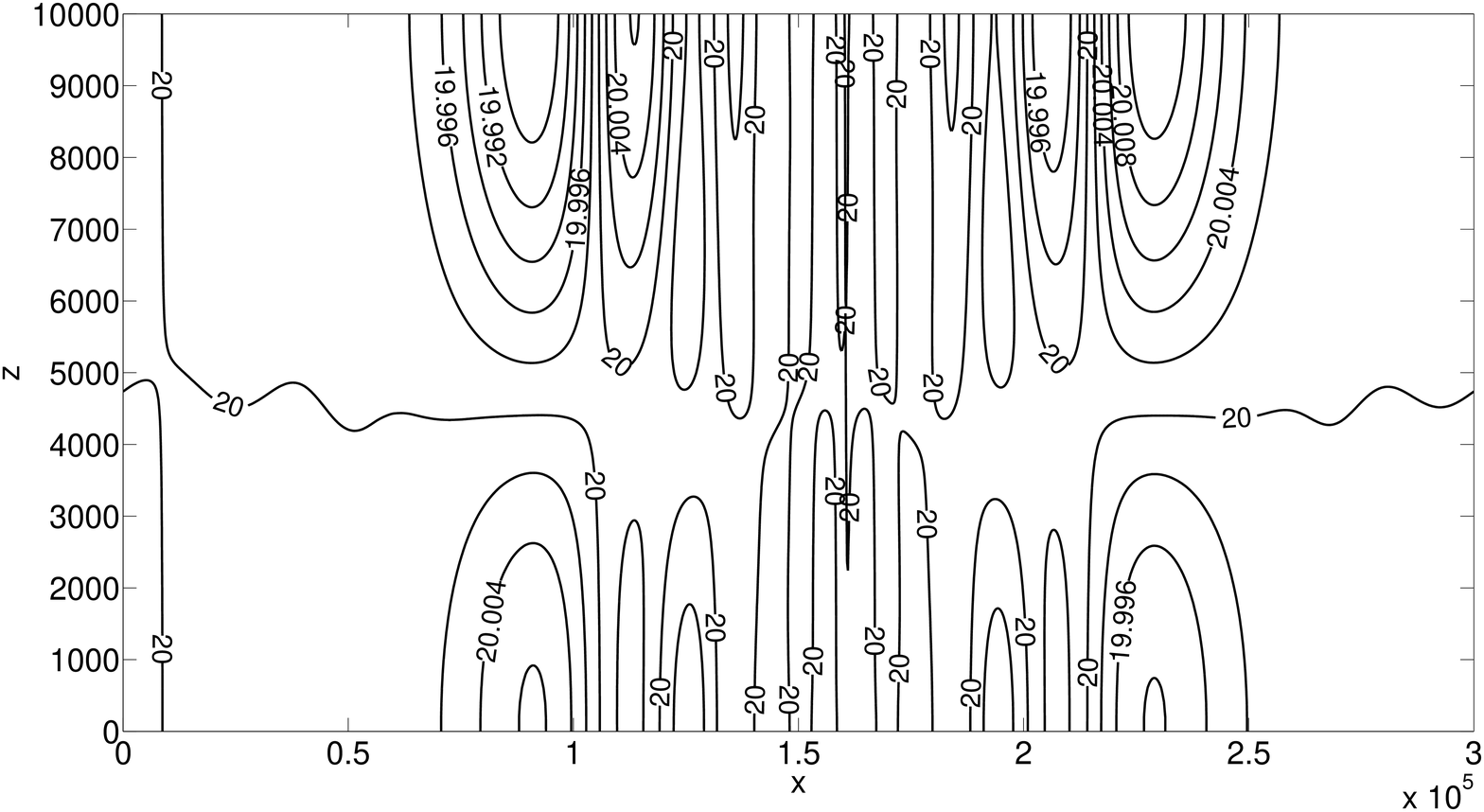}
 \caption{Contours of horizontal velocity in the internal gravity wave test.}
\label{fig:igw_u}
\end{figure}

\begin{figure}[htbc] 
 \centering
 \includegraphics[width= \textwidth]{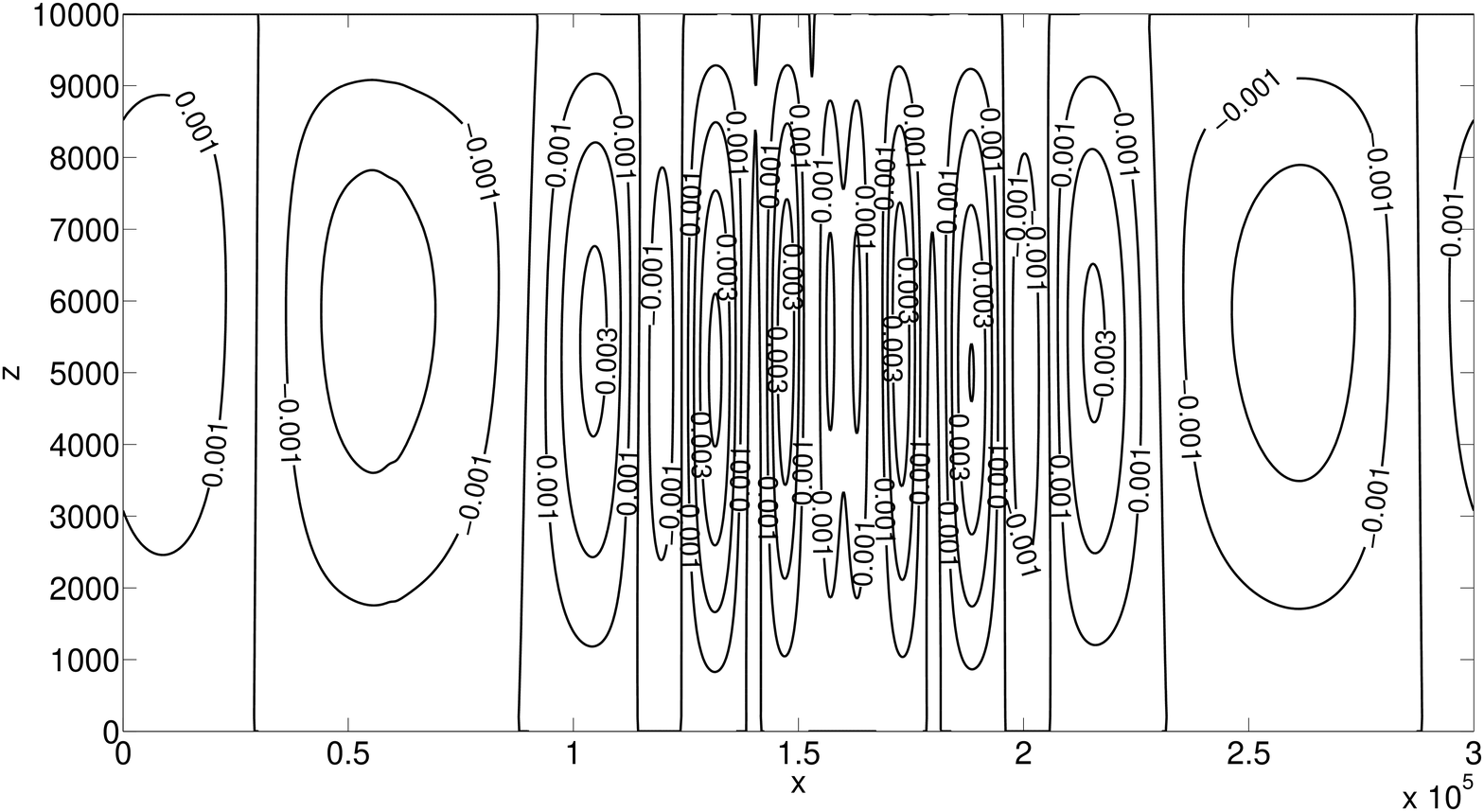}
 \caption{Contours of vertical velocity in the internal gravity wave test.}
\label{fig:igw_w}
\end{figure}

In this section we consider the test case  proposed in \cite{skamarock:1994}. 
It consists in a set of inertia-gravity waves propagating in a channel with a uniformly
stratified reference atmosphere characterized by a constant Brunt-W\"ais\"al\"a frequency
$N^2 = 0.01$. The domain and the initial and boundary conditions are identical to those of \cite{skamarock:1994}.
The initial perturbation in potential temperature radiates symmetrically to the left and to the right, but because of 
the superimposed mean horizontal flow ($u=20 $m/s), does not remain centered around the initial position.
Contours of potential temperature perturbation, horizontal velocity, and vertical velocity  time $t_f=3000 $ s are shown in figures 
\ref{fig:igw_pottemp}, \ref{fig:igw_u}, \ref{fig:igw_w}, respectively. The computed results
compare well with the structure displayed by the analytical solution of the 
linearized equations proposed in \cite{baldauf:2013} and with numerical results obtained with other numerical
methods, see e.g. \cite{bonaventura:2000}. It is to be  remarked that for this experiment 
$300\times 10$ elements, $p^{\pi}=4,$  $p^{u}=5$ and a timestep $\Delta t = 15 $ s  were used,
corresponding to a Courant number $C_{snd} \approx 25.$


\clearpage

\subsection{Rising thermal bubble}
\label{warm_bubble}
\begin{figure}[htbc]
\subfigure
  {\includegraphics[width=.55\linewidth]{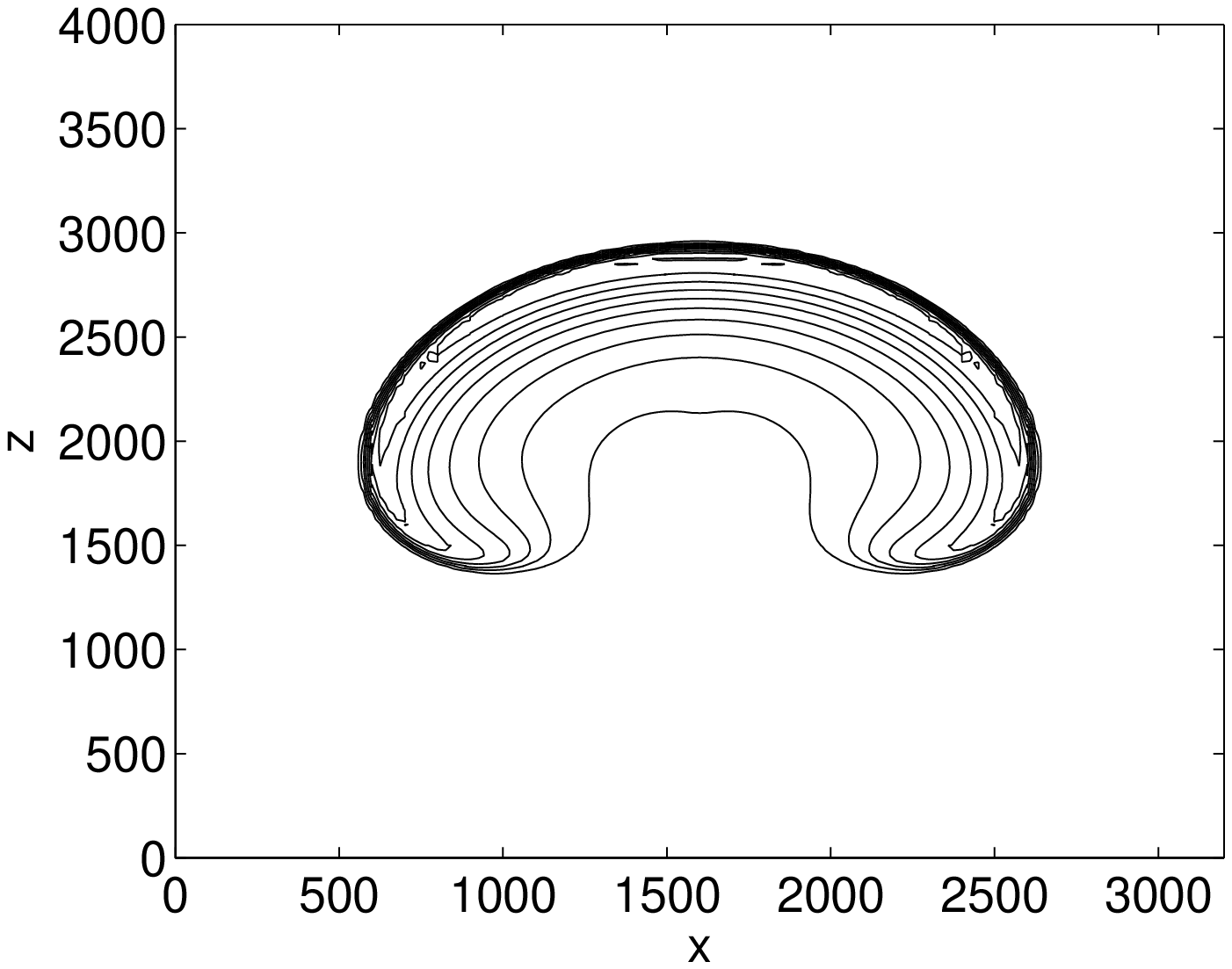}}\hfill
\subfigure
  {\includegraphics[width=.55\linewidth]{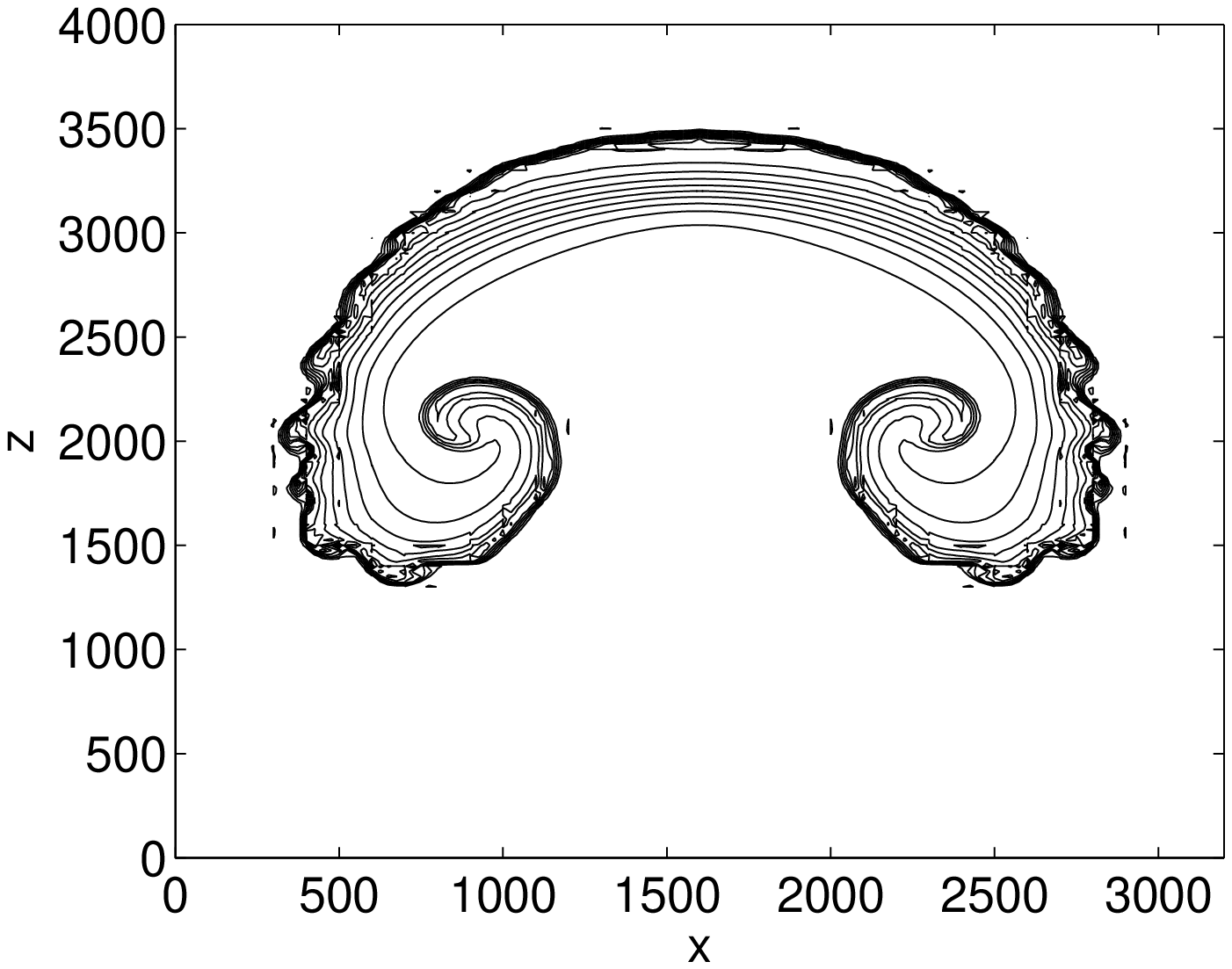}}\hfill
\\
\subfigure
  {\includegraphics[width=.55\linewidth]{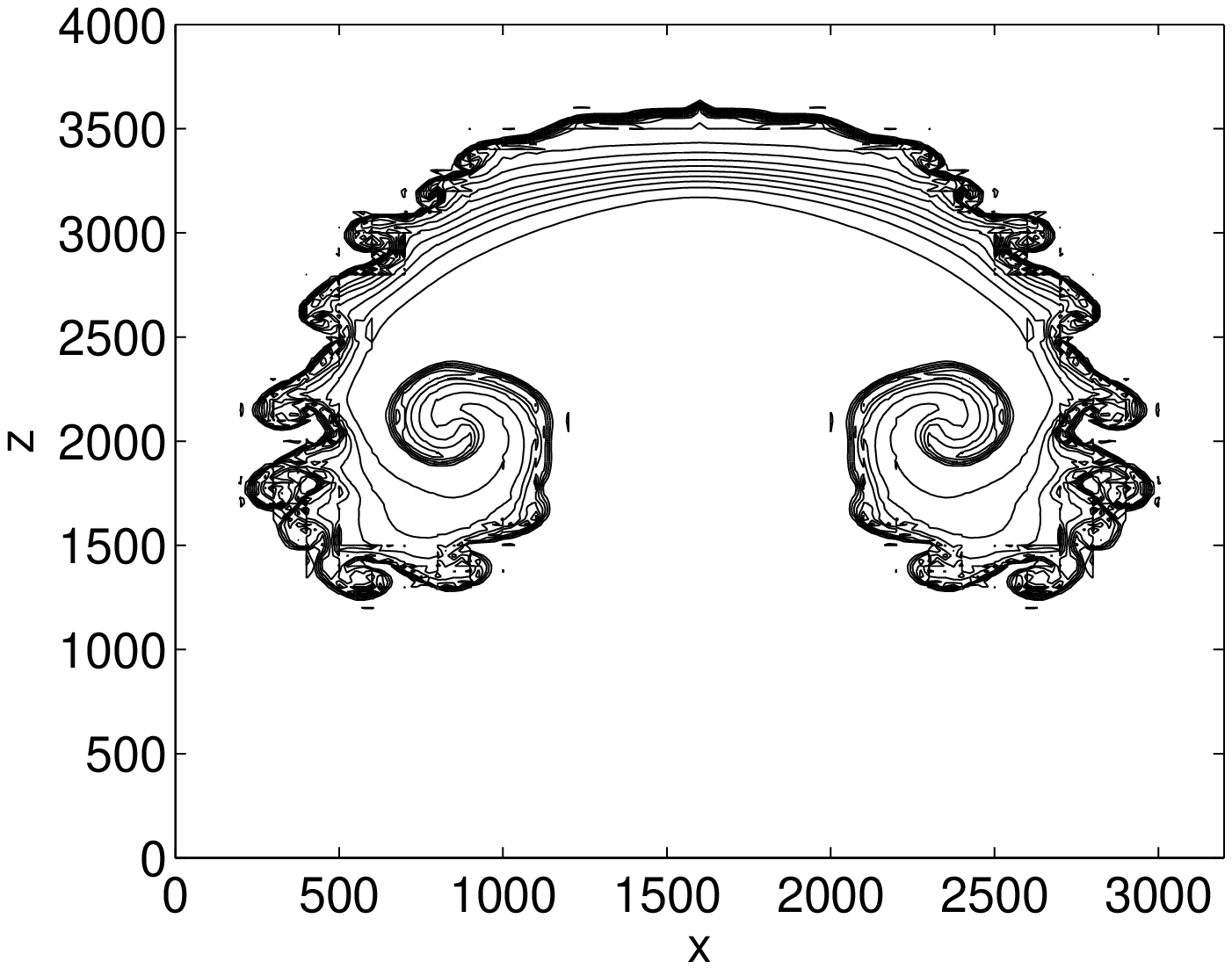}}\hfill
\subfigure
  {\includegraphics[width=.55\linewidth]{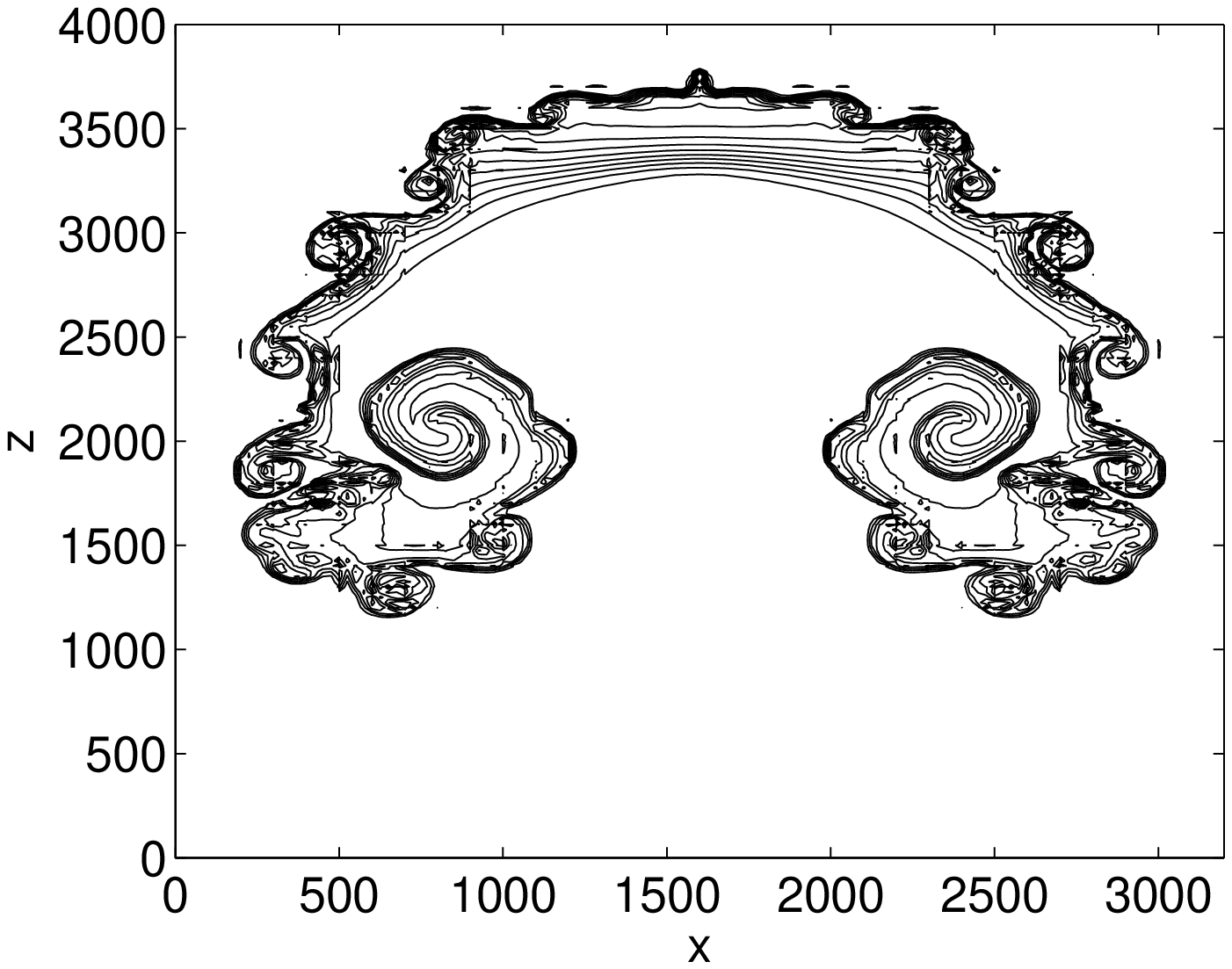}}\hfill
\caption{Contours (every 0.2 K and the zero contour is omitted) of perturbation potential temperature in the rising thermal bubble test at time 10 min, 14 min, 15 min and 16 min 
         respectively in clockwise sense.}
\label{fig:bubble_w_bn_nofill}           
\end{figure}

As nonlinear nonhydrostatic time-dependent experiment, we consider in this section the test case proposed in \cite{carpenter:1990}. 
It consists in the evolution of a warm bubble placed in an isentropic atmosphere at rest. All data are as in \cite{carpenter:1990}.
Contours of potential temperature perturbation at different times are shown in figure 
\ref{fig:bubble_w_bn_nofill}. These results were obtained using $64\times 80$ elements, $p^{\pi}=4,$  $p^{u}=5$ and a timestep $\Delta t = 0.5 $ s,
corresponding to a Courant number $C_{snd} \approx 17.$

\section{Conclusions and future perspectives}
\label{conclu}
 
 We have introduced an accurate and efficient discretization approach 
 for typical model equations of atmospheric flows. 
 We have extended to spherical geometry  the techniques proposed in \cite{tumolo:2013},
 combining  a semi-Lagrangian approach 
 with the TR-BDF2 semi-implicit time discretization method and with
 a spatial discretization based on adaptive discontinuous finite elements.
 The resulting method is unconditionally stable and has full second order accuracy 
 in time, thus improving standard off-centered trapezoidal rule discretizations
 without any major increase in the computational cost nor loss in stability,
 while allowing the use of time steps up to 100 times larger than those
 required by stability for explicit methods applied to corresponding DG discretizations.
 The method also has arbitrarily high order accuracy in space and can effectively adapt the number of degrees
 of freedom employed in each element in order to balance accuracy and computational
 cost.   The $p-$adaptivity approach employed does not require remeshing and is especially
 suitable for applications, such as numerical weather prediction, in which a large number
 of physical quantities is associated to a given the mesh.
 
 Furthermore, although the proposed method can be implemented on arbitrary unstructured and
 nonconforming meshes, like reduced Gaussian grids employed by spectral transform models,
 even in applications on simple Cartesian meshes in spherical coordinates the $p-$adaptivity approach can
 cure effectively the pole problem by reducing the polynomial degree in the polar elements, yielding a
 reduction in the computational cost that is comparable to that achieved with reduced grids.
 Numerical simulations of classical shallow water and non-hydrostatic benchmarks have been employed 
 to validate the method and to demonstrate its capability to achieve accurate results even at large Courant numbers, 
 while reducing the computational cost thanks to the adaptivity approach. 
 The proposed numerical framework can thus provide the basis of for an accurate and efficient adaptive weather prediction system.

 \pagebreak

\section*{Acknowledgements}

This research work has been  supported financially by the The Abdus Salam International Center for Theoretical Physics,
Earth System Physics Section.
We are extremely grateful to Filippo Giorgi of ICTP for his strong interest in our work and his 
continuous support. Financial support has also been provided  by  the INDAM-GNCS 2012 project
\textit{Sviluppi teorici ed applicativi dei metodi Semi-Lagrangiani} and by Politecnico di Milano.
 We would also like to acknowledge useful conversations on the topics of this paper with
C. Erath, F. X. Giraldo, M. Restelli, N. Wood.
 
\bibliographystyle{plain}
\bibliography{SISLDG}

\begin{thebibliography}{10}

\bibitem{baldauf:2013}
M.~Baldauf and S.~Brdar.
\newblock An analytic solution for linear gravity waves in a channel as a test
  for numerical models using the non-hydrostatic, compressible {E}uler
  equations.
\newblock {\em Quarterly Journal of the Royal Meteorological Society},
  139:1977--1989, 2013.

\bibitem{bank:1985}
R.E. Bank, W.M. Coughran, W.~Fichtner, E.H. Grosse, D.J. Rose, and R.K. Smith.
\newblock {T}ransient {S}imulation of {S}ilicon {D}evices and {C}ircuits.
\newblock {\em IEEE {T}ransactions on {E}lectron {D}evices.}, 32:1992--2007,
  1985.

\bibitem{bates:1990}
J.R. Bates, F.H.M. Semazzi, R.W. Higgins, and S.R.M. Barros.
\newblock {I}ntegration of the {S}hallow {W}ater {E}quations on the {S}phere
  {U}sing a {V}ector {S}emi-{L}agrangian {S}cheme with a {M}ultigrid {S}olver.
\newblock {\em Monthly Weather Review}, 118:1615--1627, 1990.

\bibitem{bonaventura:2000}
L.~Bonaventura.
\newblock A semi-implicit, semi-{L}agrangian scheme using the height coordinate
  for a nonhydrostatic and fully elastic model of atmospheric flows.
\newblock {\em Journal of Computational Physics}, 158:186--213, 2000.

\bibitem{bonaventura:2012}
L.~Bonaventura, R.~Redler, and R.~Budich.
\newblock {\em Earth System Modelling 2: Algorithms, Code Infrastructure and
  Optimisation.}
\newblock Springer Verlag, New York, 2012.

\bibitem{butcher:2000}
J.~Butcher and D.~Chen.
\newblock A new type of singly-implicit {R}unge-{K}utta method.
\newblock {\em Applied numerical mathematics}, 34:179--188, 2000.

\bibitem{casulli:1994}
V.~Casulli and E.~Cattani.
\newblock Stability, accuracy and efficiency of a semi-implicit method for
  three-dimensional shallow water flow.
\newblock {\em Computational Mathematics and Applications}, 27:99--112, 1994.

\bibitem{cote:1988}
J.~Cot\'e.
\newblock A {L}agrange multiplier approach for the metric terms of
  semi-{L}agrangian models on the sphere.
\newblock {\em Quarterly Journal of the Royal Meteorological Society},
  114:1347--1352, 1988.

\bibitem{cote:1988b}
J.~Cot\'e and A.~Staniforth.
\newblock A {T}wo-{T}ime-{L}evel {S}emi-{L}agrangian {S}emi-implicit {S}cheme
  for {S}pectral {M}odels.
\newblock {\em Monthly Weather Review}, 116:2003--2012, 1988.

\bibitem{cullen:1990}
M.J.P. Cullen.
\newblock A test of a semi-implicit integration technique for a fully
  compressible non-hydrostatic model.
\newblock {\em Quarterly Journal of the Royal Meteorological Society},
  116:1253--1258, 1990.

\bibitem{davies:2005}
T.~Davies, M.J.P. Cullen, A.J. Malcolm, M.H. Mawson, A.~Staniforth, A.A. White,
  and N.~Wood.
\newblock A new dynamical core for the {M}et {O}ffice's global and regional
  modelling of the atmosphere.
\newblock {\em Quarterly Journal of the Royal Meteorological Society},
  131:1759--1782, 2005.

\bibitem{dawson:2006}
C.N. Dawson, J.J. Westerink, J.C. Feyen, and D.~Pothina.
\newblock Continuous, {D}iscontinuous and coupled {D}iscontinuous-{C}ontinuous
  {G}alerkin finite element methods for the shallow water equations.
\newblock {\em International Journal of Numerical Methods in Fluids},
  52:63--88, 2006.

\bibitem{desharnais:1990}
F.~Desharnais and A.~Robert.
\newblock Errors near the poles generated by a semi-{L}agrangian integration
  scheme in a global spectral model.
\newblock {\em Atmosphere-Ocean}, 28:162--176, 1990.

\bibitem{dumbser:2013}
M.~Dumbser and V.~Casulli.
\newblock A staggered semi-implicit spectral discontinuous {G}alerkin scheme
  for the shallow water equations.
\newblock {\em Applied Mathematics and Computation}, 219(15):8057 -- 8077,
  2013.

\bibitem{gill:1982}
A.E. Gill.
\newblock {\em {A}tmospheric-{O}cean {D}ynamics.}
\newblock {A}cademic {P}ress, 1987.

\bibitem{giraldo:1999}
F.X. Giraldo.
\newblock Trajectory computations for spherical geodesic grids in cartesian
  space.
\newblock {\em Monthly Weather Review}, 127:1651--1662, 1999.

\bibitem{giraldo:2002}
F.X. Giraldo, J.S. Hesthaven, and T.~Warburton.
\newblock High-{O}rder {D}iscontinuous {G}alerkin methods for the spherical
  shallow water equations.
\newblock {\em Journal of Computational Physics}, 181:499--525, 2002.

\bibitem{giraldo:2013}
F.X. Giraldo, J.F. Kelly, and E.M. Constantinescu.
\newblock Implicit-explicit formulations of a three-dimensional nonhydrostatic
  unified model of the atmosphere ({N}{U}{M}{A}).
\newblock {\em SIAM Journal of Scientific Computing}, 35(5):1162--1194, 2013.

\bibitem{giraldo:2010}
F.X. Giraldo and M.~Restelli.
\newblock High-order semi-implicit time-integrators for a triangular
  discontinuous {G}alerkin oceanic shallow water model.
\newblock {\em International Journal of Numerical Methods in Fluids},
  63(9):1077--1102, 2010.

\bibitem{hortal:1991}
M.~Hortal and A.~Simmons.
\newblock Use of reduced {G}aussian grids in spectral models.
\newblock {\em Monthly Weather Review}, 119:1057--1074, 1991.

\bibitem{hosea:1996}
M.E. Hosea and L.F. Shampine.
\newblock Analysis and implementation of {TR}-{BDF}2.
\newblock {\em Applied Numerical Mathematics}, 20:21--37, 1996.

\bibitem{jakob:1995}
R.~Jakob-Chien, J.~Hack, and D.~Williamson.
\newblock Spectral transform solutions to the shallow water test set.
\newblock {\em Journal of Computational Physics}, 119:164Ð187, 1995.

\bibitem{carpenter:1990}
R.~L.~Carpenter Jr., K.~K. Droegemeier, P.~R. Woodward, and C.~E. Hane.
\newblock {A}pplication of the {P}iecewise {P}arabolic {M}ethod ({PPM}) to
  {M}eteorological {M}odeling.
\newblock {\em Monthly Weather Review}, 118:586--612, 1990.

\bibitem{kelley:1995}
C.~T. Kelley.
\newblock {\em Iterative {M}ethods for {L}inear and {N}onlinear {E}quations.}
\newblock SIAM, Philadelphia, 1995.

\bibitem{kelly:2012}
J.F. Kelly and F.X. Giraldo.
\newblock Continuous and discontinuous {G}alerkin methods for a scalable
  three-dimensional nonhydrostatic atmospheric model: Limited-area mode.
\newblock {\em Journal of Computational Physics}, 231(24):7988--8008, 2012.

\bibitem{kennedy:2003}
C.~Kennedy and M.~Carpenter.
\newblock Additive {R}unge-{K}utta schemes for convection-diffusion-reaction
  equations.
\newblock {\em Applied Numerical Mathematics}, 44:139--181, 2003.

\bibitem{lambert:1991}
J.D. Lambert.
\newblock {\em Numerical methods for ordinary differential systems}.
\newblock Wiley, 1991.

\bibitem{lauter:2008}
M.~L{\"a}uter, F.X. Giraldo, D.~Handorf, and K.~Dethloff.
\newblock A discontinuous {G}alerkin method for the shallow water equations in
  spherical triangular coordinates.
\newblock {\em Journal of Computational Physics}, 227(24):10226--10242, 2008.

\bibitem{lauter:2005}
M.~L\"auter, D.~Handorf, and K.~Dethloff.
\newblock Unsteady analytical solutions of the spherical shallow water
  equations.
\newblock {\em Journal of Computational Physics}, 210:535,553, 2005.

\bibitem{leroux:2013}
D.Y. {Le Roux}.
\newblock Spurious inertial oscillations in shallow-water models.
\newblock {\em Journal of Computational Physics}, 231:7959--7987, 2013.

\bibitem{leroux:2005}
D.Y. {Le Roux} and G.F. Carey.
\newblock Stability-dispersion analysis of the discontinuous {G}alerkin
  linearized shallow-water system.
\newblock {\em International Journal of Numerical Methods in Fluids},
  48:325--347, 2005.

\bibitem{leveque:2007}
R.~LeVeque.
\newblock {\em Finite difference methods for ordinary and partial differential
  equations: steady-state and time-dependent problems.}
\newblock Society for Industrial and Applied Mathematics, 2007.

\bibitem{mcdonald:1989}
A.~McDonald and J.R. Bates.
\newblock Semi-{L}agrangian {I}ntegration of a {G}ridpoint {S}hallow {W}ater
  {M}odel on the {S}phere.
\newblock {\em Monthly Weather Review}, 117:130--137, 1989.

\bibitem{mcgregor:1993}
J.L. McGregor.
\newblock Economical determination of departure points for semi-{L}agrangian
  models.
\newblock {\em Monthly Weather Review}, 121:221--330, 1993.

\bibitem{morton:1998}
K.~W. Morton.
\newblock On the analysis of finite volume methods for evolutionary problems.
\newblock {\em SIAM Journal of Numerical Analysis}, 35:2195--2222, 1998.

\bibitem{morton:1988}
K.~W. Morton, A.~Priestley, and E.~S\"uli.
\newblock Stability of the {L}agrange-{G}alerkin scheme with inexact
  integration.
\newblock {\em RAIRO Modellisation Matemathique et Analyse Numerique},
  22:625--653, 1988.

\bibitem{morton:1995}
K.~W. Morton and E.~S\"uli.
\newblock Evolution-{G}alerkin methods and their supraconvergence.
\newblock {\em Numerische Mathematik}, 71:331--355, 1995.

\bibitem{nair:2005b}
R.~D. Nair, S.J. Thomas, and R.D. Loft.
\newblock A {D}iscontinuous {G}alerkin {G}lobal {S}hallow {W}ater {M}odel.
\newblock {\em Monthly Weather Review}, 133:876--888, 2005.

\bibitem{nair:2005}
R.~D. Nair, S.J. Thomas, and R.D. Loft.
\newblock A {D}iscontinuous {G}alerkin transport scheme on the cubed sphere.
\newblock {\em Monthly Weather Review}, 133:814--828, 2005.

\bibitem{priestley:1994}
A.~Priestley.
\newblock {E}xact {P}rojections and the {L}agrange-{G}alerkin {M}ethod: {A}
  {R}ealistic {A}lternative to {Q}uadrature.
\newblock {\em Journal of Computational Physics}, 112:316--333, 1994.

\bibitem{restelli:2006}
M.~Restelli, L.~Bonaventura, and R.~Sacco.
\newblock A semi-{L}agrangian {D}iscontinuous {G}alerkin method for scalar
  advection by incompressible flows.
\newblock {\em Journal of Computational Physics}, 216:195--215, 2006.

\bibitem{restelli:2009}
M.~Restelli and F.X. Giraldo.
\newblock A conservative {D}iscontinuous {G}alerkin semi-implicit formulation
  for the {N}avier-{S}tokes equations in nonhydrostatic mesoscale modeling.
\newblock {\em SIAM Journal of Scientific Computing}, 31:2231--2257, 2009.

\bibitem{ripodas:2009}
P.~Ripodas, A.~Gassmann, J.~F\"orstner, D.~Majewski, M.~Giorgetta, P.~Korn,
  L.~Kornblueh, H.~Wan, G.~Z\"angl, L.~Bonaventura, and T.~Heinze.
\newblock Icosahedral {S}hallow {W}ater {M}odel ({ICOSWM}): results of shallow
  water test cases and sensitivity to model parameters.
\newblock {\em Geoscientific Model Development}, 2:231--251, 2009.

\bibitem{ritchie:1988}
H.~Ritchie.
\newblock Application of the {S}emi-{L}agrangian {M}ethod to a {S}pectral
  {M}odel of the {S}hallow {W}ater {E}quations.
\newblock {\em Monthly Weather Review}, 116:1587--1598, 1988.

\bibitem{rosatti:2005}
G.~Rosatti, L.~Bonaventura, and D.~Cesari.
\newblock Semi-implicit, semi-{L}agrangian environmental modelling on cartesian
  grids with cut cells.
\newblock {\em Journal of Computational Physics}, 204:353--377, 2005.

\bibitem{saad:1986}
Y.~Saad and M.H. Schultz.
\newblock {GMRES}: A generalized minimal residual algorithm for solving
  nonsymmetric linear systems.
\newblock {\em SIAM Journal on Scientific and Statistical Computing},
  7:856--869, 1986.

\bibitem{skamarock:1994}
W.C. Skamarock and J.B. Klemp.
\newblock {E}fficiency and accuracy of the {K}lemp-{W}ilhelmson time-splitting
  technique.
\newblock {\em Monthly Weather Review}, 122:2623--2630, 1994.

\bibitem{staniforth:2010}
A.~Staniforth, A.A. White, and N.~Wood.
\newblock Treatment of vector equations in deep-atmosphere semi-{L}agrangian
  models.{I}: {M}omentum equation.
\newblock {\em Quarterly Journal of the Royal Meteorological Society},
  136:497--506, 2010.

\bibitem{temperton:2001}
C.~Temperton, M.~Hortal, and A.~Simmons.
\newblock A two-time-level semi-{L}agrangian global spectral model.
\newblock {\em Quarterly Journal of the Royal Meteorological Society},
  127:111--127, 2001.

\bibitem{thuburn:2000}
J.~Thuburn and Y.~Li.
\newblock Numerical simulations of {R}ossby-{H}aurwitz waves.
\newblock {\em Tellus A}, 52:181–189, 2000.

\bibitem{thuburn:2013}
J.~Thuburn and A.A. White.
\newblock A geometrical view of the shallow-atmosphere approximation, with
  application to the semi-lagrangian departure point calculation.
\newblock {\em Quarterly Journal of the Royal Meteorological Society},
  139:261–268, 2013.

\bibitem{tumolo:2013}
G.~Tumolo, L.~Bonaventura, and M.~Restelli.
\newblock A semi-implicit, semi-{L}agrangian, p-adaptive discontinuous
  {G}alerkin method for the shallow water equations.
\newblock {\em Journal of Computational Physics}, 232:46--67, January 2013.

\bibitem{walters:1983a}
R.A. Walters.
\newblock Numerically induced oscillations in finite-element approximations to
  the shallow-water equations.
\newblock {\em International Journal of Numerical Methods in Fluids},
  3:591--604, 1983.

\bibitem{walters:1983b}
R.A. Walters and G.F. Carey.
\newblock Analysis of spurious oscillation modes for the shallow-water and
  {N}avier-{S}tokes equations.
\newblock {\em Computers and {F}luids}, 11:51--68, 1983.

\bibitem{williamson:1992}
D.L. Williamson, J.B. Drake, J.J. Hack, R.~Jacob, and P.N. Swarztrauber.
\newblock {A} {S}tandard {T}est {S}et for the {N}umerical {A}pproximations to
  the {S}hallow {W}ater {E}quations in {S}pherical {G}eometry.
\newblock {\em Journal of Computational Physics}, 102:211--224, 1992.

\bibitem{zienkiewicz:1983}
O.C. Zienkiewicz, J.P.Gago, and D.W. Kelly.
\newblock The hierarchical concept in finite element analysis.
\newblock {\em Computers and {S}tructures}, 16:53--65, 1983.

\end{thebibliography}

\end{document}